\documentclass[12pt]{article}
\usepackage{amssymb}
\usepackage{amsmath, amsthm}
\usepackage{epsfig}
\usepackage{graphicx}

\makeatletter
\@addtoreset{equation}{section}
\makeatother

\newtheorem{theorem}{Theorem}[section]

\newtheorem{conjecture}[theorem]{Conjecture}
\theoremstyle{definition}

\newtheorem{remark}[theorem]{Remark}

\begin{document}

\title{Upper bounds on the smallest size of a complete arc in the plane $%
PG(2,q)$}
\date{}
\author{Daniele Bartoli \\
{\footnotesize Dipartimento di Matematica e Informatica, Universit\`{a}
degli Studi di Perugia, }\\
{\footnotesize Via Vanvitelli~1, Perugia, 06123, Italy}\\
{\footnotesize e-mail: daniele.bartoli@dmi.unipg.it}\and
Alexander A.
Davydov \\
{\footnotesize Institute for Information Transmission Problems, Russian
Academy of Sciences, }\\
{\footnotesize Bol'shoi Karetnyi per. 19, GSP-4, Moscow, 127994, Russia}\\
{\footnotesize e-mail: adav@iitp.ru}
\and Giorgio Faina , Stefano Marcugini, Fernanda Pambianco \\
{\footnotesize Dipartimento di Matematica e Informatica, Universit\`{a}
degli Studi di Perugia, }\\
{\footnotesize Via Vanvitelli~1, Perugia, 06123,
Italy}\\
{\footnotesize e-mail: faina@dmi.unipg.it (G. Faina),
gino@dmi.unipg.it (S. Marcugini)},\\{\footnotesize
fernanda@dmi.unipg.it (F. Pambianco)}} \maketitle

\begin{abstract}
New upper bounds on the smallest size $t_{2}(2,q)$ of a
complete arc in the projective plane $PG(2,q)$ are obtained for
$q\leq 9109.$ From these new bounds it follows that for $q\leq
2621$ and $q=2659,2663,2683,2693,2753,2801$, the relation
$t_{2}(2,q)<4.5\sqrt{q}$ holds. Also, for $q\leq 5399$ and
$q=5413,5417,5419,5441,5443,5471,5483,5501,5521,$ we have $
t_{2}(2,q)<4.8\sqrt{q}.$ Finally, for $q\leq 9067$ it holds
that $ t_{2}(2,q)<5\sqrt{q}.$ The new upper bounds are obtained
by finding new small complete arcs with the help of a computer
search using randomized greedy algorithms.
\end{abstract}


\section{Introduction}

Let $PG(2,q)$ be the projective plane over the Galois field
$F_{q}$. An $n$ -arc is a set of $n$ points no three of which
are collinear. An $n$-arc is called complete if it is not
contained in an $(n+1)$-arc of $PG(2,q)$. For an introduction
in projective geometries over finite fields, see \cite
{HirsBook},\cite{SegreLeGeom},\cite{SegreIntrodGalGeom}.

In \cite{HirsSt-old},\cite{HirsStor-2001} the close relationship between the
theory of $n$-arcs, coding theory and mathematical statistics is presented.
In particular, a complete arc in a plane $PG(2,q),$ points of which are
treated as 3-dimensional $q$-ary columns, defines a parity check matrix of a
$q$-ary linear code with codimension 3, Hamming distance 4, and covering
radius 2. Arcs can be interpreted as linear maximum distance separable (MDS)
codes \cite{szoT93},\cite{thaJ92d} and they are related to optimal coverings
arrays \cite{Hartman-Haskin} and to superregular matrices \cite{Keri}.

One of the main problems in the study of projective planes, which is also of
interest in coding theory, is finding of the spectrum of possible sizes of
complete arcs. In particular, the value of $t_{2}(2,q)$, the smallest size
of a complete arc in $PG(2,q),$ is interesting. Finding an estimation of the
minimum size $t_{2}(2,q)$ is a hard open problem.

This work is devoted to \emph{upper bounds }on $t_{2}(2,q)$.

Surveys of results on the sizes of plane complete arcs, methods
of their construction and comprehension of the relating
properties can be found in
\cite{BDFMP-DM}-\cite{DFMP-JG2009},\cite{FP},\cite{HirsSurvey83}-\cite
{HirsStor-2001},\cite{SegreLeGeom}-\cite{szoT93}.

Problems connected with small complete plane arcs are
considered in \cite
{BDFMP-DM}-\cite{DFMP-JG2009},\cite{DGMP-Innov},\cite{DGMP-JCD},\cite
{FainaGiul},
\cite{Giul2000}-\cite{GiulUghi},\cite{KV},\cite{LisMarcPamb2008}
,\cite{MMP-q29}-\cite{Ost},\cite{Polv},\cite{SegreLeGeom},\cite{SZ},
see also the references therein.

The exact values of $t_{2}(2,q)$ are known only for $q\leq 32,$
see \cite {MMP-q29} and recent work \cite{MaMiP32} where the
equalities $ t_{2}(2,31)=t_{2}(2,32)=14$ are proven. Also,
there are the following lower bounds (see
\cite{Polv},\cite{SegreLeGeom}):
\begin{equation*}
t_{2}(2,q)>\left\{
\begin{array}{ll}
\sqrt{2q}+1 & \text{for any }q \\
\sqrt{3q}+\frac{1}{2} & \text{for }q=p^{h},\text{ }p\text{ prime, }h=1,2,3
\end{array}
\right. .
\end{equation*}

Let $t(\mathcal{P}_{q})$ be the size of the smallest complete arc in any
(not necessarily Galois) projective plane $\mathcal{P}_{q}$ of order $q$. In
\cite{KV}, for \emph{sufficiently large} $q$, the following result is proved
(we give it in the form of \cite[Tab.\thinspace 2.6]{HirsStor-2001}):
\begin{equation}
t(\mathcal{P}_{q})\leq d\sqrt{q}\log ^{c}q,\text{ }c\leq 300,
\label{eq1_KimVu_c=300}
\end{equation}
where $c$ and $d$ are constants independent of $q$ (i.e. universal
constants). The logarithm basis is not noted as the estimate is asymptotic.

Following to \cite{BDFMP-DM}, we denote the aggregates of $q$ values:
\begin{eqnarray*}
T_{2}
&=&\{5119,5147,5153,5209,5231,5237,5261,5279,5281,5303,5347,5641,5843,6011,
\\
&&8192\}; \\
T_{3} &=&\{2^{14},2^{15},2^{18}\}; \\
Q &=&\{961,1024,1369,1681,2401\}=\{31^{2},2^{10},37^{2},41^{2},7^{4}\}.
\end{eqnarray*}

Let $\overline{t}_{2}(2,q)$ be the smallest \emph{known }size of a complete
arc in $PG(2,q)$.

For $q\leq 841$ the values of $\overline{t}_{2}(2,q)$ (up to
June 2009) are collected in \cite[ Tab.\thinspace
1]{DFMP-JG2009} whence it follows that $
\overline{t}_{2}(2,q)<4\sqrt{q}$ for $q\leq 841$. In
\cite{Giul2000}, see also \cite{GiulUghi}, complete
$(4\sqrt{q}-4)$-arcs are obtained for $ q=p^{2} $ odd, $q\leq
1681$ or $q=2401$. In \cite{DGMP-JCD},\cite{FainaGiul} complete
$(4\sqrt{q}-4)$-arcs are obtained for even $q=64,256,1024.$

In this work we showed that $\overline{t}_{2}(2,857)=117<4\sqrt{857}.$

So, it holds that
\begin{equation}
t_{2}(2,q)<4\sqrt{q}\quad \text{for }2\leq q\leq 841,\text{ }q=857,\text{ }
q\in Q.  \label{eq1_<4sqroot(q)}
\end{equation}

For $q\leq 5107$ and $q\in T_{2}\cup T_{3}$ the values of
$\overline{t} _{2}(2,q)$ (up to June 2011) are collected in
\cite[Tabs.\thinspace 1-4] {BDFMP-DM} where the following
results are obtained:
\begin{eqnarray}
t_{2}(2,q) &<&4.5\sqrt{q}~~\,\mbox{ for }q\leq 2593,\,\,q=2693,2753.
\label{eq1_4.5} \\
t_{2}(2,q) &<&4.79\sqrt{q}\,\,\mbox{ for }q\leq 5107.  \label{eq1_4.79} \\
t_{2}(2,q) &<&4.98\sqrt{q}\,\,\mbox{ for }q\in T_{2}.  \label{eq1_4.98}
\end{eqnarray}
\begin{equation}
t_{2}(2,q)<0.9987\sqrt{q}\ln ^{0.75}q\quad \mbox{ for }23\leq q\leq 5107,
\text{ }q\in T_{2}\cup T_{3}.  \label{eq1_9987}
\end{equation}

Moreover, in \cite{BDFMP-DM} the following conjectures were proposed:

\begin{conjecture}
\label{conj1_ ln0.75}\emph{\ \cite{BDFMP-DM} }In $PG(2,q),$
\begin{equation}
t_{2}(2,q)<\sqrt{q}\ln ^{0.75}q\quad \mbox{ for }q\geq 23.
\label{eq1_conj0.75}
\end{equation}
\end{conjecture}

\begin{conjecture}
\label{conj1_5sqrt_q} \emph{\ \cite{BDFMP-DM} }In $PG(2,q),$
\begin{equation}
t_{2}(2,q)<5\sqrt{q}\hspace{1cm}\quad \mbox{ for }q\leq 8192.
\label{eq1_conj5sqrt_q}
\end{equation}
\end{conjecture}

In this work we obtained many new small arcs and extended and improved
results of \cite{BDFMP-DM}. We have proved Conjecture \ref{conj1_5sqrt_q},
see Theorem \ref{th1_<4.5_5}.

Results of this work allow us to hope that Conjectures
\ref{conj1_ ln0.75} is true.

We denote the aggregates of $q$ values:
\begin{eqnarray*}
T_{4} &=&\{359,367,401,419,512,541,571,643,653,719,773,787\}; \\
T_{5}
&=&
\{857,881,919,929,941,953,967,1019,1031,1069,1097,1109,1123,1151,1163,1187,
\\
&&1201,1217,1231,1259,1289,1301,1319,1331,1361,1373,1433,1447,1493,1511,1523,
\\
&&1553,1567,1571,1583,1597,1601,1613,1627,1663,1693,1697,1723,1741,1759,1777,
\\
&&1789,1823,1871,1873,1889,1907,1973,1987,1993,2003,2039,2111,2113,2129,2131,
\\
&&2141,2143,2179,2197,2213,2237,2251,2269,2287,2309,2339,2341,2357,2399,2411,
\\
&&2417,2437,2467,2473,2531,2609,2617,2621\}; \\
T_{6}
&=&
\{2657,2659,2663,2677,2683,2699,2719,2741,2797,2801,2819,2833,2837,2851,2857,
\\
&&2879,2897,2917,2953,2957,2971,2999,3011,3019,3037,3041,3061,3137,3181,3217,
\\
&&3221,3259,3307,3329,3331,3371,3373,3391,3407,3449,3461,3527,3541,3547,3557,
\\
&&3581,3613,3631,3671,3673,3677,3691,3697,3701,3719,3721,3761,3767,3823,3833,
\\
&&3847,3851,3877,3917,3923,3943,3947,3989,4007,4051,4079,4096,4127,4129,4201,
\\
&&4337,4339,4391,4409,4451,4483,4507,4603,4621,4673,4729,4751,4793,4799,4903,
\\
&&4931,4973,4999,5023,5051,5077,5081,5099,5101,5153,5209,5231,5261,5279,5281,
\\
&&5347,5641,6011,8192\}.
\end{eqnarray*}

In this paper we obtained complete arcs with sizes smaller than
in \cite {BDFMP-DM} (i.e. we improved upper bounds on
$t_{2}(2,q)$) for $q\in T_{4}$ (in the region $q\leq 841),$ for
$q\in T_{5}$ (in the region $853\leq q\leq 2621)$, and for
$q\in T_{6}$ (in the region $2633\leq q)$.

\begin{theorem}
\label{th1_<4.5_5} In $PG(2,q),$ the following holds.
\begin{eqnarray*}
t_{2}(2,q) &<&4.5\sqrt{q}~~\,\mbox{ for }q\leq
2621,\,\,q=2659,2663,2683,2693,2753,2801. \\
t_{2}(2,q) &<&4.8\sqrt{q}~~\,\mbox{ for }q\leq
5399,\,\,q=5413,5417,5419,5441,5443,5471,5483,5501,5521. \\
t_{2}(2,q) &<&5\sqrt{q}\,\,\quad ~\mbox{ for }q\leq 9067.
\end{eqnarray*}
\end{theorem}

\begin{theorem}
\label{th1_ln0.75} In $PG(2,q),$
\begin{equation}
t_{2}(2,q)<0.9987\sqrt{q}\ln ^{0.75}q\quad \mbox{ for }23\leq q\leq 9109,
\text{ }q\in T_{3}.  \label{eq1_ln0.75}
\end{equation}
\end{theorem}

In whole, this work can be considered as development of work
\cite{BDFMP-DM}.

\section{Small complete $k$-arcs in $PG(2,q)$, $q\leq 9109$}

Throughout the paper, in all tables we denote
$A_{q}=\left\lfloor a_{q}\sqrt{
q}-\overline{t}_{2}(2,q)\right\rfloor $ where
\begin{equation*}
a_{q}=\left\{
\begin{array}{cl}
4 & \text{if }q\leq 841,\text{ }q\in Q \\
4.5 & \text{if }853\leq q\leq
2621,\,\,q=2659,2663,2683,2693,2753,2801,\,q\notin Q \\
5 & \text{if }2623\leq q\leq 9067,\,q\notin
\{2659,2663,2683,2693,2753,2801\}.
\end{array}
\right. .
\end{equation*}
Also, in all tables, $B_{q}$ is a superior approximation of
$\overline{t} _{2}(2,q)/\sqrt{q}.$

For $q\leq 841$, the values of $\overline{t}_{2}(2,q)$ (up to
June 2011) are collected in \cite[Tab.\thinspace 1]{BDFMP-DM}.
In this work we obtained small arcs with new sizes for $q\in
T_{1}.$ The new arcs are obtained by computer search, based on
the randomized greedy algorithms. The current values of
$\overline{t}_{2}(2,q)$ for $q\leq 841$ are given in Table 1.
The data for $q\in T_{1}$ improving results of
\cite[Tab.\thinspace 1]{BDFMP-DM} are written in Table~1 in
bold font. The exact values $\overline{t} _{2}(2,q)=t_{2}(2,q)$
are marked by the dot~\textquotedblleft $\centerdot $
\textquotedblright . In particular, due to the recent result
\cite{MaMiP32} we noted the values
$t_{2}(2,31)=t_{2}(2,32)=14.$

From Table 1 and the results of \cite{Giul2000,GiulUghi}, on
complete $(4 \sqrt{q}-4)$-arcs for $q=p^{2}$ (see Introduction)
we obtain Theorem \ref {Th3_4sqroot(q)} improving and extending
the results of \cite[Th.\thinspace 3.1]{BDFMP-DM}.

\newpage

\textbf{Table 1. }The smallest known sizes
$\overline{t}_{2}=\overline{t} _{2}(2,q)<4\sqrt{q}$ of complete
arcs in planes $\mathrm{PG}(2,q),$ $q\leq 841.$
$A_{q}=\left\lfloor
4\sqrt{q}-\overline{t}_{2}(2,q)\right\rfloor $, $ B_{q}\geq
\overline{t}_{2}(2,q)/\sqrt{q}$\smallskip

$\renewcommand{\arraystretch}{0.88}
\begin{array}{@{}r@{\,\,\,\,\,}l@{\,\,\,\,}r@{\,\,\,\,}c|c@{\,\,\,\,}c@{\,\,\,\,}r@{\,\,\,\,}c|cc@{\,\,\,\,}cc|cr@{\,\,\,\,}rc}
\hline
q & \overline{t}_{2} & A_{q} & B_{q} & q & \overline{t}_{2} & A_{q} & B_{q}
& q & \overline{t}_{2} & A_{q} & B_{q} & q & \overline{t}_{2} & A_{q} &
B_{q}^{\phantom{H^{L}}} \\ \hline
2 & \phantom{1}4\centerdot & 1 & 2.83 & 128 & 34 & 11 & 3.01 & 347 & 67 & 7
& 3.60 & 599 & 94 & 3 & 3.85 \\
3 & \phantom{1}4\centerdot & 2 & 2.31 & 131 & 36 & 9 & 3.15 & 349 & 67 & 7 &
3.59 & 601 & 90 & 8 & 3.68 \\
4 & \phantom{1}6\centerdot & 2 & 3.00 & 137 & 37 & 9 & 3.17 & 353 & 68 & 7 &
3.62 & 607 & 95 & 3 & 3.86 \\
5 & \phantom{1}6\centerdot & 2 & 2.69 & 139 & 37 & 10 & 3.14 & \mathbf{359}
& \mathbf{68} & \mathbf{7} & \mathbf{3.59} & 613 & 96 & 3 & 3.88 \\
7 & \phantom{1}6\centerdot & 4 & 2.27 & 149 & 39 & 9 & 3.20 & 361 & 69 & 7 &
3.64 & 617 & 96 & 3 & 3.87 \\
8 & \phantom{1}6\centerdot & 5 & 2.13 & 151 & 39 & 10 & 3.18 & \mathbf{367}
& \mathbf{69} & \mathbf{7} & \mathbf{3.61} & 619 & 96 & 3 & 3.86 \\
9 & \phantom{1}6\centerdot & 6 & 2.00 & 157 & 40 & 10 & 3.20 & 373 & 70 & 7
& 3.63 & 625 & 96 & 4 & 3.84 \\
11 & \phantom{1}7\centerdot & 6 & 2.12 & 163 & 41 & 10 & 3.22 & 379 & 71 & 6
& 3.65 & 631 & 97 & 3 & 3.87 \\
13 & \phantom{1}8\centerdot & 6 & 2.22 & 167 & 42 & 9 & 3.26 & 383 & 71 & 7
& 3.63 & 641 & 98 & 3 & 3.88 \\
16 & \phantom{1}9\centerdot & 7 & 2.25 & 169 & 42 & 10 & 3.24 & 389 & 72 & 6
& 3.66 & \mathbf{643} & \mathbf{98} & \mathbf{3} & \mathbf{3.87} \\
17 & 10\centerdot & 6 & 2.43 & 173 & 43 & 9 & 3.27 & 397 & 73 & 6 & 3.67 &
647 & 99 & 2 & 3.90 \\
19 & 10\centerdot & 7 & 2.30 & 179 & 44 & 9 & 3.29 & \mathbf{401} & \mathbf{
73} & \mathbf{7} & \mathbf{3.65} & \mathbf{653} & \mathbf{99} & \mathbf{3} &
\mathbf{3.88} \\
23 & 10\centerdot & 9 & 2.09 & 181 & 44 & 9 & 3.28 & 409 & 74 & 6 & 3.66 &
659 & 100 & 2 & 3.90 \\
25 & 12\centerdot & 8 & 2.40 & 191 & 46 & 9 & 3.33 & \mathbf{419} & \mathbf{
75} & \mathbf{6} & \mathbf{3.67} & 661 & 90 & 12 & 3.51 \\
27 & 12\centerdot & 8 & 2.31 & 193 & 46 & 9 & 3.32 & 421 & 76 & 6 & 3.71 &
673 & 101 & 2 & 3.90 \\
29 & 13\centerdot & 8 & 2.42 & 197 & 47 & 9 & 3.35 & 431 & 77 & 6 & 3.71 &
677 & 102 & 2 & 3.93 \\
31 & 14\centerdot & 8 & 2.52 & 199 & 47 & 9 & 3.34 & 433 & 77 & 6 & 3.71 &
683 & 102 & 2 & 3.91 \\
32 & 14\centerdot & 8 & 2.48 & 211 & 49 & 9 & 3.38 & 439 & 78 & 5 & 3.73 &
691 & 103 & 2 & 3.92 \\
37 & 15 & 9 & 2.47 & 223 & 51 & 8 & 3.42 & 443 & 78 & 6 & 3.71 & 701 & 104 &
1 & 3.93 \\
41 & 16 & 9 & 2.50 & 227 & 51 & 9 & 3.39 & 449 & 79 & 5 & 3.73 & 709 & 104 &
2 & 3.91 \\
43 & 16 & 10 & 2.45 & 229 & 51 & 9 & 3.38 & 457 & 80 & 5 & 3.75 & \mathbf{719
} & \mathbf{105} & \mathbf{2} & \mathbf{3.92} \\
47 & 18 & 9 & 2.63 & 233 & 52 & 9 & 3.41 & 461 & 80 & 5 & 3.73 & 727 & 106 &
1 & 3.94 \\
49 & 18 & 10 & 2.58 & 239 & 53 & 8 & 3.43 & 463 & 80 & 6 & 3.72 & 729 & 104
& 4 & 3.86 \\
53 & 18 & 11 & 2.48 & 241 & 53 & 9 & 3.42 & 467 & 81 & 5 & 3.75 & 733 & 107
& 1 & 3.96 \\
59 & 20 & 10 & 2.61 & 243 & 53 & 9 & 3.40 & 479 & 82 & 5 & 3.75 & 739 & 107
& 1 & 3.94 \\
61 & 20 & 11 & 2.57 & 251 & 55 & 8 & 3.48 & 487 & 83 & 5 & 3.77 & 743 & 108
& 1 & 3.97 \\
64 & 22 & 10 & 2.75 & 256 & 55 & 9 & 3.44 & 491 & 83 & 5 & 3.75 & 751 & 108
& 1 & 3.95 \\
67 & 23 & 9 & 2.81 & 257 & 55 & 9 & 3.44 & 499 & 84 & 5 & 3.77 & 757 & 109 &
1 & 3.97 \\
71 & 22 & 11 & 2.62 & 263 & 56 & 8 & 3.46 & 503 & 85 & 4 & 3.79 & 761 & 109
& 1 & 3.96 \\
73 & 24 & 10 & 2.81 & 269 & 57 & 8 & 3.48 & 509 & 85 & 5 & 3.77 & 769 & 110
& 0 & 3.97 \\
79 & 26 & 9 & 2.93 & 271 & 57 & 8 & 3.47 & \mathbf{512} & \mathbf{85} &
\mathbf{5} & \mathbf{3.76} & \mathbf{773} & \mathbf{110} & \mathbf{1} &
\mathbf{3.96} \\
81 & 26 & 10 & 2.89 & 277 & 58 & 8 & 3.49 & 521 & 86 & 5 & 3.77 & \mathbf{787
} & \mathbf{111} & \mathbf{1} & \mathbf{3.96} \\
83 & 27 & 9 & 2.97 & 281 & 59 & 8 & 3.52 & 523 & 86 & 5 & 3.77 & 797 & 112 &
0 & 3.97 \\
89 & 28 & 9 & 2.97 & 283 & 59 & 8 & 3.51 & 529 & 87 & 5 & 3.79 & 809 & 113 &
0 & 3.98 \\
97 & 30 & 9 & 3.05 & 289 & 60 & 8 & 3.53 & \mathbf{541} & \mathbf{88} &
\mathbf{5} & \mathbf{3.79} & 811 & 113 & 0 & 3.97 \\
101 & 30 & 10 & 2.99 & 293 & 60 & 8 & 3.51 & 547 & 89 & 4 & 3.81 & 821 & 114
& 0 & 3.98 \\
103 & 31 & 9 & 3.06 & 307 & 62 & 8 & 3.54 & 557 & 90 & 4 & 3.82 & 823 & 114
& 0 & 3.98 \\
107 & 32 & 9 & 3.10 & 311 & 63 & 7 & 3.58 & 563 & 91 & 3 & 3.84 & 827 & 115
& 0 & 4.00 \\
109 & 32 & 9 & 3.07 & 313 & 63 & 7 & 3.57 & 569 & 91 & 4 & 3.82 & 829 & 115
& 0 & 4.00 \\
113 & 33 & 9 & 3.11 & 317 & 63 & 8 & 3.54 & \mathbf{571} & \mathbf{91} &
\mathbf{4} & \mathbf{3.81} & 839 & 115 & 0 & 3.98 \\
121 & 34 & 10 & 3.10 & 331 & 65 & 7 & 3.58 & 577 & 92 & 4 & 3.84 & 841 & 112
& 4 & 3.87 \\
125 & 35 & 9 & 3.14 & 337 & 66 & 7 & 3.60 & 587 & 93 & 3 & 3.84 &  &  &  &
\\
127 & 35 & 10 & 3.11 & 343 & 66 & 8 & 3.57 & 593 & 94 & 3 & 3.87 &  &  &  &
\\ \hline
\end{array}
$

\begin{theorem}
\label{Th3_4sqroot(q)} In $PG(2,q),$ the following holds.
\begin{eqnarray}
t_{2}(2,q) &<&\phantom{.5}4\sqrt{q}~\mbox{ for }2\leq q\leq 841,\text{ }
q=857,\text{ }q\in Q.  \label{eq3_<4sqroot(q)} \\
t_{2}(2,q) &\leq &\phantom{.5}3\sqrt{q}~\mbox{ for }2\leq q\leq 89,\text{ }
\,~q=101;  \notag \\
t_{2}(2,q) &<&3.5\sqrt{q}~\mbox{ for }2\leq q\leq 277;  \notag \\
t_{2}(2,q) &<&3.6\sqrt{q}~\mbox{ for }2\leq q\leq 349,\text{ }q=359,661;
\notag \\
t_{2}(2,q) &<&3.7\sqrt{q}~\mbox{ for }2\leq q\leq 419,\text{ }q=601,661;
\notag \\
t_{2}(2,q) &<&3.8\sqrt{q}~\mbox{ for }2\leq q\leq 541,~q=601,661;  \notag \\
t_{2}(2,q) &<&3.9\sqrt{q}~\mbox{ for }2\leq q\leq 673,~q=729,961,1024.
\notag
\end{eqnarray}
Also,
\begin{eqnarray*}
t_{2}(2,q) &\leq &4\sqrt{q}-9\ \mbox{ for }\,\,37\leq q\leq 211,\text{ }
q=23,227,229,233,241,243,256,257,661; \\
t_{2}(2,q) &\leq &4\sqrt{q}-8\ \mbox{ for }\,\,23\leq q\leq 307,\text{ }
q=317,343,601,661; \\
t_{2}(2,q) &\leq &4\sqrt{q}-7\ \mbox{ for }\ 19\leq q\leq 373,\
q=383,401,601,661; \\
t_{2}(2,q) &\leq &4\sqrt{q}-6\ \mbox{ for
}\,\,\phantom{1}9\leq q\leq 433,\ q=443,463,601,661; \\
t_{2}(2,q) &\leq &4\sqrt{q}-5\ \mbox{ for }\ \phantom{1}8\leq q\leq 499,\
q=509,512,521,523,529,541,601,661; \\
t_{2}(2,q) &\leq &4\sqrt{q}-4\ \mbox{ for }\ \phantom{1}7\leq q\leq 557,\
q=569,571,577,601,625,661,729,841,\,q\in Q; \\
t_{2}(2,q) &<&4\sqrt{q}-3\ \mbox{ for }\ \phantom{1}7\leq q\leq 643,\
q=653,661,729,841,\,q\in Q; \\
t_{2}(2,q) &\leq &4\sqrt{q}-2\ \mbox{ for }\ \phantom{1}3\leq q\leq 691,
\text{ }q=709,719,729,841,\,q\in Q; \\
t_{2}(2,q) &<&4\sqrt{q}-1\ \mbox{ for }\ \phantom{1}2\leq q\leq 761,\text{ }
q=773,787,841,\,q\in Q.
\end{eqnarray*}
\end{theorem}

For $853\leq q\leq 2621$, the values of $\overline{t}_{2}(2,q)$
(up to June 2011) are collected in \cite[Tabs\thinspace
2,3]{BDFMP-DM}. In this work we obtained small arcs with new
sizes for $q\in T_{5}.$ The new arcs are obtained by computer
search, based on the randomized greedy algorithms. The current
values of $\overline{t}_{2}(2,q)<4.5\sqrt{q}$\, for $853\leq
q\leq 2621 $ are given in Table 2. The data for $q\in T_{5}$
improving results of \cite[Tabs\thinspace 2,3]{BDFMP-DM} are
written in Table~2 in bold font. The data for $q=p^{2}$ with
$\overline{t}_{2}(2,q)=4\sqrt{q}-4$ \cite {Giul2000,GiulUghi}
and for $q=857$ with $\overline{t}_{2}(2,857)=117<4\sqrt{ 857}$
are written in italic font.

From Table 2 and the results of \cite{Giul2000,GiulUghi}, on
complete $(4 \sqrt{q}-4)$-arcs for $q=p^{2}$ (see Introduction)
we obtain Theorem \ref {Th3_4.5sqroot(q)} improving and
extending the results of \cite[ Th.\thinspace
3.2]{BDFMP-DM}.\newpage

\begin{center}
\textbf{Table 2.} The smallest known sizes
$\overline{t}_{2}=\overline{t} _{2}(2,q)<4.5\sqrt{q}$ of
complete arcs in planes $PG(2,q),$ $853\leq q\leq 2621,$
$A_{q}=\left\lfloor
a_{q}\sqrt{q}-\overline{t}_{2}(2,q)\right\rfloor $ , $B_{q}\geq
\overline{t}_{2}(2,q)/\sqrt{q}$\smallskip

$\renewcommand{\arraystretch}{0.9}
\begin{array}{@{}r@{\,\,\,}c@{\,\,\,}c@{\,\,\,\,}c|@{\,\,}c@{\,\,\,}c@{\,\,\,\,}c@{\,\,\,}c|@{\,\,}c@{\,\,\,\,}c@{\,\,\,\,}c@{\,\,\,}c|@{\,\,}c@{\,\,\,\,}c@{\,\,\,\,}c@{\,\,\,\,}c}
\hline
q & \overline{t}_{2} & A_{q} & B_{q} & q & \overline{t}_{2} & A_{q} & B_{q}
& q & \overline{t}_{2} & A_{q} & B_{q} & q & \overline{t}_{2} & A_{q} &
B_{q}^{\phantom{H^{L}}} \\ \hline
853 & 117 & 14 & 4.01 & 1279 & 150 & 10 & 4.20 & 1699 & 178 & 7 & 4.32 & 2161
& 205 & 4 & 4.41 \\
\mathbf{857} & \emph{117} & \emph{0} & \emph{4.00} & 1283 & 150 & 11 & 4.19
& 1709 & 178 & 8 & 4.31 & \mathbf{2179} & \mathbf{206} & \mathbf{4} &
\mathbf{4.42} \\
859 & 118 & 13 & 4.03 & \mathbf{1289} & \mathbf{150} & \mathbf{11} & \mathbf{
4.18} & 1721 & 179 & 7 & 4.32 & 2187 & 207 & 3 & 4.43 \\
863 & 118 & 14 & 4.02 & 1291 & 151 & 10 & 4.21 & \mathbf{1723} & \mathbf{179}
& \mathbf{7} & \mathbf{4.32} & \mathbf{2197} & \mathbf{207} & \mathbf{3} &
\mathbf{4.42} \\
877 & 119 & 14 & 4.02 & 1297 & 151 & 11 & 4.20 & 1733 & 180 & 7 & 4.33 & 2203
& 208 & 3 & 4.44 \\
\mathbf{881} & \mathbf{119} & \mathbf{14} & \mathbf{4.01} & \mathbf{1301} &
\mathbf{151} & \mathbf{11} & \mathbf{4.19} & \mathbf{1741} & \mathbf{180} &
\mathbf{7} & \mathbf{4.32} & 2207 & 208 & 3 & 4.43 \\
883 & 120 & 13 & 4.04 & 1303 & 151 & 11 & 4.19 & 1747 & 181 & 7 & 4.34 & 2209
& 208 & 3 & 4.43 \\
887 & 120 & 14 & 4.03 & 1307 & 152 & 10 & 4.21 & 1753 & 181 & 7 & 4.33 &
\mathbf{2213} & \mathbf{208} & \mathbf{3} & \mathbf{4.43} \\
907 & 122 & 13 & 4.06 & \mathbf{1319} & \mathbf{152} & \mathbf{11} & \mathbf{
4.19} & \mathbf{1759} & \mathbf{181} & \mathbf{7} & \mathbf{4.32} & 2221 &
209 & 3 & 4.44 \\
911 & 122 & 13 & 4.05 & 1321 & 153 & 10 & 4.21 & \mathbf{1777} & \mathbf{182}
& \mathbf{7} & \mathbf{4.32} & \mathbf{2237} & \mathbf{209} & \mathbf{3} &
\mathbf{4.42} \\
\mathbf{919} & \mathbf{122} & \mathbf{14} & \mathbf{4.03} & 1327 & 153 & 10
& 4.21 & 1783 & 183 & 7 & 4.34 & 2239 & 210 & 2 & 4.44 \\
\mathbf{929} & \mathbf{123} & \mathbf{14} & \mathbf{4.04} & \mathbf{1331} &
\mathbf{153} & \mathbf{11} & \mathbf{4.20} & 1787 & 183 & 7 & 4.33 & 2243 &
210 & 3 & 4.44 \\
937 & 124 & 13 & 4.06 & \mathbf{1361} & \mathbf{155} & \mathbf{11} & \mathbf{
4.21} & \mathbf{1789} & \mathbf{183} & \mathbf{7} & \mathbf{4.33} & \mathbf{
2251} & \mathbf{210} & \mathbf{3} & \mathbf{4.43} \\
\mathbf{941} & \mathbf{124} & \mathbf{14} & \mathbf{4.05} & 1367 & 156 & 10
& 4.22 & 1801 & 184 & 6 & 4.34 & 2267 & 211 & 3 & 4.44 \\
947 & 125 & 13 & 4.07 & \emph{1369} & \emph{144} & \emph{4} & \emph{3.90} &
1811 & 184 & 7 & 4.33 & \mathbf{2269} & \mathbf{211} & \mathbf{3} & \mathbf{
4.43} \\
\mathbf{953} & \mathbf{125} & \mathbf{13} & \mathbf{4.05} & \mathbf{1373} &
\mathbf{156} & \mathbf{10} & \mathbf{4.22} & \mathbf{1823} & \mathbf{185} &
\mathbf{7} & \mathbf{4.34} & 2273 & 212 & 2 & 4.45 \\
\emph{961} & \emph{120} & \emph{4} & \emph{3.88} & 1381 & 157 & 10 & 4.23 &
1831 & 186 & 6 & 4.35 & 2281 & 212 & 2 & 4.44 \\
\mathbf{967} & \mathbf{126} & \mathbf{13} & \mathbf{4.06} & 1399 & 158 & 10
& 4.23 & 1847 & 187 & 6 & 4.36 & \mathbf{2287} & \mathbf{212} & \mathbf{3} &
\mathbf{4.44} \\
971 & 127 & 13 & 4.08 & 1409 & 159 & 9 & 4.24 & 1849 & 187 & 6 & 4.35 & 2293
& 213 & 2 & 4.45 \\
977 & 127 & 13 & 4.07 & 1423 & 160 & 9 & 4.25 & 1861 & 188 & 6 & 4.36 & 2297
& 213 & 2 & 4.45 \\
983 & 128 & 13 & 4.09 & 1427 & 160 & 9 & 4.24 & 1867 & 188 & 6 & 4.36 &
\mathbf{2309} & \mathbf{213} & \mathbf{3} & \mathbf{4.44} \\
991 & 127 & 14 & 4.04 & 1429 & 160 & 10 & 4.24 & \mathbf{1871} & \mathbf{188}
& \mathbf{6} & \mathbf{4.35} & 2311 & 214 & 2 & 4.46 \\
997 & 129 & 13 & 4.09 & \mathbf{1433} & \mathbf{160} & \mathbf{10} & \mathbf{
4.23} & \mathbf{1873} & \mathbf{188} & \mathbf{6} & \mathbf{4.35} & 2333 &
215 & 2 & 4.46 \\
1009 & 130 & 12 & 4.10 & 1439 & 161 & 9 & 4.25 & 1877 & 189 & 5 & 4.37 &
\mathbf{2339} & \mathbf{215} & \mathbf{2} & \mathbf{4.45} \\
1013 & 130 & 13 & 4.09 & \mathbf{1447} & \mathbf{161} & \mathbf{10} &
\mathbf{4.24} & 1879 & 189 & 6 & 4.37 & \mathbf{2341} & \mathbf{214} &
\mathbf{3} & \mathbf{4.43} \\
\mathbf{1019} & \mathbf{130} & \mathbf{13} & \mathbf{4.08} & 1451 & 162 & 9
& 4.26 & \mathbf{1889} & \mathbf{189} & \mathbf{6} & \mathbf{4.35} & 2347 &
216 & 2 & 4.46 \\
1021 & 131 & 12 & 4.10 & 1453 & 162 & 9 & 4.25 & 1901 & 190 & 6 & 4.36 & 2351
& 216 & 2 & 4.46 \\
\emph{1024} & \emph{124} & \emph{4} & \emph{3.88} & 1459 & 162 & 9 & 4.25 &
\mathbf{1907} & \mathbf{190} & \mathbf{6} & \mathbf{4.36} & \mathbf{2357} &
\mathbf{216} & \mathbf{2} & \mathbf{4.45} \\
\mathbf{1031} & \mathbf{131} & \mathbf{13} & \mathbf{4.08} & 1471 & 163 & 9
& 4.25 & 1913 & 191 & 5 & 4.37 & 2371 & 217 & 2 & 4.46 \\
1033 & 132 & 12 & 4.11 & 1481 & 164 & 9 & 4.27 & 1931 & 192 & 5 & 4.37 & 2377
& 216 & 3 & 4.44 \\
1039 & 132 & 13 & 4.10 & 1483 & 164 & 9 & 4.26 & 1933 & 192 & 5 & 4.37 & 2381
& 217 & 2 & 4.45 \\
1049 & 133 & 12 & 4.11 & 1487 & 164 & 9 & 4.26 & 1949 & 193 & 5 & 4.38 & 2383
& 218 & 1 & 4.47 \\
1051 & 133 & 12 & 4.11 & 1489 & 164 & 9 & 4.26 & 1951 & 193 & 5 & 4.37 & 2389
& 218 & 1 & 4.47 \\
1061 & 134 & 12 & 4.12 & \mathbf{1493} & \mathbf{164} & \mathbf{9} & \mathbf{
4.25} & \mathbf{1973} & \mathbf{194} & \mathbf{5} & \mathbf{4.37} & 2393 &
218 & 2 & 4.46 \\
1063 & 134 & 12 & 4.11 & 1499 & 165 & 9 & 4.27 & 1979 & 195 & 5 & 4.39 &
\mathbf{2399} & \mathbf{218} & \mathbf{2} & \mathbf{4.46} \\
\mathbf{1069} & \mathbf{134} & \mathbf{13} & \mathbf{4.10} & \mathbf{1511} &
\mathbf{165} & \mathbf{9} & \mathbf{4.25} & \mathbf{1987} & \mathbf{195} &
\mathbf{5} & \mathbf{4.38} & \emph{2401} & \emph{192} & \emph{4} & \emph{3.92
} \\
1087 & 136 & 12 & 4.13 & \mathbf{1523} & \mathbf{166} & \mathbf{9} & \mathbf{
4.26} & \mathbf{1993} & \mathbf{195} & \mathbf{5} & \mathbf{4.37} & \mathbf{
2411} & \mathbf{219} & \mathbf{1} & \mathbf{4.47} \\
1091 & 136 & 12 & 4.12 & 1531 & 167 & 9 & 4.27 & 1997 & 196 & 5 & 4.39 &
\mathbf{2417} & \mathbf{219} & \mathbf{2} & \mathbf{4.46} \\
1093 & 136 & 12 & 4.12 & 1543 & 167 & 9 & 4.26 & 1999 & 196 & 5 & 4.39 & 2423
& 220 & 1 & 4.47 \\
\mathbf{1097} & \mathbf{136} & \mathbf{13} & \mathbf{4.11} & 1549 & 168 & 9
& 4.27 & \mathbf{2003} & \mathbf{196} & \mathbf{5} & \mathbf{4.38} & \mathbf{
2437} & \mathbf{220} & \mathbf{2} & \mathbf{4.46} \\
1103 & 137 & 12 & 4.13 & \mathbf{1553} & \mathbf{168} & \mathbf{9} & \mathbf{
4.27} & 2011 & 197 & 4 & 4.40 & 2441 & 221 & 1 & 4.48 \\ \hline
\end{array}
$
\end{center}

\newpage

\begin{center}
\textbf{Table 2 }(continue)\medskip

$\renewcommand{\arraystretch}{1.0}
\begin{array}{@{}r@{\,\,\,\,}c@{\,\,\,\,}c@{\,\,\,\,}c|@{\,\,}c@{\,\,\,\,}c@{\,\,\,\,}c@{\,\,\,}c|@{\,\,}c@{\,\,\,\,}c@{\,\,\,\,}c@{\,\,\,}c|@{\,\,}c@{\,\,\,\,}c@{\,\,\,\,}c@{\,\,\,\,}c}
\hline
q & \overline{t}_{2} & A_{q} & B_{q} & q & \overline{t}_{2} & A_{q} & B_{q}
& q & \overline{t}_{2} & A_{q} & B_{q} & q & \overline{t}_{2} & A_{q} &
B_{q}^{\phantom{H^{L}}} \\ \hline
\mathbf{1109} & \mathbf{137} & \mathbf{12} & \mathbf{4.12} & 1559 & 169 & 8
& 4.29 & 2017 & 197 & 5 & 4.39 & 2447 & 221 & 1 & 4.47 \\
1117 & 138 & 12 & 4.13 & \mathbf{1567} & \mathbf{169} & \mathbf{9} & \mathbf{
4.27} & 2027 & 198 & 4 & 4.40 & 2459 & 222 & 1 & 4.48 \\
\mathbf{1123} & \mathbf{138} & \mathbf{12} & \mathbf{4.12} & \mathbf{1571} &
\mathbf{169} & \mathbf{9} & \mathbf{4.27} & 2029 & 198 & 4 & 4.40 & \mathbf{
2467} & \mathbf{222} & \mathbf{1} & \mathbf{4.47} \\
1129 & 139 & 12 & 4.14 & 1579 & 170 & 8 & 4.28 & \mathbf{2039} & \mathbf{198}
& \mathbf{5} & \mathbf{4.39} & \mathbf{2473} & \mathbf{222} & \mathbf{1} &
\mathbf{4.47} \\
\mathbf{1151} & \mathbf{140} & \mathbf{12} & \mathbf{4.13} & \mathbf{1583} &
\mathbf{170} & \mathbf{9} & \mathbf{4.28} & 2048 & 199 & 4 & 4.40 & 2477 &
223 & 0 & 4.49 \\
1153 & 141 & 11 & 4.16 & \mathbf{1597} & \mathbf{171} & \mathbf{8} & \mathbf{
4.28} & 2053 & 199 & 4 & 4.40 & 2503 & 224 & 1 & 4.48 \\
\mathbf{1163} & \mathbf{141} & \mathbf{12} & \mathbf{4.14} & \mathbf{1601} &
\mathbf{171} & \mathbf{9} & \mathbf{4.28} & 2063 & 200 & 4 & 4.41 & 2521 &
225 & 0 & 4.49 \\
1171 & 142 & 11 & 4.15 & 1607 & 172 & 8 & 4.30 & 2069 & 200 & 4 & 4.40 &
\mathbf{2531} & \mathbf{225} & \mathbf{1} & \mathbf{4.48} \\
1181 & 143 & 11 & 4.17 & 1609 & 172 & 8 & 4.29 & 2081 & 201 & 4 & 4.41 & 2539
& 226 & 0 & 4.49 \\
\mathbf{1187} & \mathbf{143} & \mathbf{12} & \mathbf{4.16} & \mathbf{1613} &
\mathbf{172} & \mathbf{8} & \mathbf{4.29} & 2083 & 201 & 4 & 4.41 & 2543 &
226 & 0 & 4.49 \\
1193 & 144 & 11 & 4.17 & 1619 & 173 & 8 & 4.30 & 2087 & 201 & 4 & 4.40 & 2549
& 226 & 1 & 4.48 \\
\mathbf{1201} & \mathbf{144} & \mathbf{11} & \mathbf{4.16} & 1621 & 173 & 8
& 4.30 & 2089 & 201 & 4 & 4.40 & 2551 & 227 & 0 & 4.50 \\
1213 & 145 & 11 & 4.17 & \mathbf{1627} & \mathbf{173} & \mathbf{8} & \mathbf{
4.29} & 2099 & 202 & 4 & 4.41 & 2557 & 227 & 0 & 4.49 \\
\mathbf{1217} & \mathbf{145} & \mathbf{11} & \mathbf{4.16} & 1637 & 174 & 8
& 4.31 & \mathbf{2111} & \mathbf{202} & \mathbf{4} & \mathbf{4.40} & 2579 &
228 & 0 & 4.49 \\
1223 & 146 & 11 & 4.18 & 1657 & 175 & 8 & 4.30 & \mathbf{2113} & \mathbf{202}
& \mathbf{4} & \mathbf{4.40} & 2591 & 229 & 0 & 4.50 \\
1229 & 146 & 11 & 4.17 & \mathbf{1663} & \mathbf{175} & \mathbf{8} & \mathbf{
4.30} & \mathbf{2129} & \mathbf{203} & \mathbf{4} & \mathbf{4.40} & 2593 &
229 & 0 & 4.50 \\
\mathbf{1231} & \mathbf{146} & \mathbf{11} & \mathbf{4.17} & 1667 & 176 & 7
& 4.32 & \mathbf{2131} & \mathbf{203} & \mathbf{4} & \mathbf{4.40} & \mathbf{
2609} & \mathbf{229} & \mathbf{0} & \mathbf{4.49} \\
1237 & 147 & 11 & 4.18 & 1669 & 176 & 7 & 4.31 & 2137 & 204 & 4 & 4.42 &
\mathbf{2617} & \mathbf{230} & \mathbf{0} & \mathbf{4.50} \\
1249 & 148 & 11 & 4.19 & \emph{1681} & \emph{160} & \emph{4} & \emph{3.91} &
\mathbf{2141} & \mathbf{204} & \mathbf{4} & \mathbf{4.41} & \mathbf{2621} &
\mathbf{230} & \mathbf{0} & \mathbf{4.50} \\
\mathbf{1259} & \mathbf{148} & \mathbf{11} & \mathbf{4.18} & \mathbf{1693} &
\mathbf{177} & \mathbf{8} & \mathbf{4.31} & \mathbf{2143} & \mathbf{204} &
\mathbf{4} & \mathbf{4.41} &  &  &  &  \\
1277 & 150 & 10 & 4.20 & \mathbf{1697} & \mathbf{177} & \mathbf{8} & \mathbf{
4.30} & 2153 & 205 & 3 & 4.42 &  &  &  &  \\ \hline
\end{array}
$\newpage
\end{center}

\begin{theorem}
\label{Th3_4.5sqroot(q)} In $PG(2,q),$ the following holds.
\begin{eqnarray}
t_{2}(2,q) &<&4.5\sqrt{q}~\mbox{ for }q\leq
2621,\,q=2659,2663,2683,2693,2753,2801.  \label{eq3_<4.5sqroot(q)} \\
t_{2}(2,q) &<&4.1\sqrt{q}~\mbox{ for }q\leq
1031,\,q=1039,1069,1369,1681,2401;  \notag \\
t_{2}(2,q) &<&4.2\sqrt{q}~\mbox{ for }q\leq
1289,\,q=1297,1301,1303,1319,1331,1369,1681,2401;  \notag \\
t_{2}(2,q) &<&4.3\sqrt{q}~\mbox{ for }q\leq
1627,\,q=1657,1663,1681,1697,2401;  \notag
\label{eq3_4.5sqroot(q)_a*sqroot(q)} \\
t_{2}(2,q) &<&4.4\sqrt{q}~\mbox{ for
}q\leq 2053,\,q=2069,2087,2089,2111,2113,2129,2131,2401.  \notag
\end{eqnarray}
Also,
\begin{eqnarray*}
t_{2}(2,q) &<&4.5\sqrt{q}-13\,\mbox{ for }q\leq \text{ }997,\text{\thinspace
}q=1013,1019,1024,1031,1039,1069,1097,1369,1681, \\
&&\phantom{4.5\sqrt{q}-9~\mbox{ for }q\leq 3911,\,q=\,\,\;}2401; \\
t_{2}(2,q) &<&4.5\sqrt{q}-12\,\mbox{ for }q\leq 1151,\text{\thinspace }
q=1163,1187,1369,1681,2401; \\
t_{2}(2,q) &<&4.5\sqrt{q}-11\,\mbox{ for }q\leq
1259,\,q=1283,1289,1297,1301,1303,1331,1319,1331,1361, \\
&&\phantom{4.5\sqrt{q}-9~\mbox{ for }q\leq 3911,\,q=\,\,\;}1369,1681,2401; \\
t_{2}(2,q) &<&4.5\sqrt{q}-10\,\mbox{ for }q\leq
1399,\,q=1429,1433,1447,1681,2401; \\
t_{2}(2,q) &<&4.5\sqrt{q}-9\,~\,\mbox{ for }q\leq
1553,\,q=1567,1571,1583,1601,1681,2401; \\
t_{2}(2,q) &<&4.5\sqrt{q}-8\,~\,\mbox{ for }q\leq
1663,\,q=1681,1693,1697,1709,2401; \\
t_{2}(2,q) &<&4.5\sqrt{q}-7\,~\,\mbox{ for }q\leq 1789,\,q=1811,1823,2401; \\
t_{2}(2,q) &<&4.5\sqrt{q}-6\,~\,\mbox{ for
}q\leq 1873,\,q=1879,1889,1901,1907,2401; \\
t_{2}(2,q) &<&4.5\sqrt{q}-5\,~\,\mbox{ for
}q\leq 2003,\,q=2017,2039,2401; \\
t_{2}(2,q) &<&4.5\sqrt{q}-4\,~\,\mbox{ for
}q\leq 2143,\,q=2161,2179,2401; \\
t_{2}(2,q) &<&4.5\sqrt{q}-3\,~\,\mbox{ for }q\leq
2237,\,q=2243,2251,2267,2269,2287,2309,2341,2377,2401; \\
t_{2}(2,q) &<&4.5\sqrt{q}-2\,~\,\mbox{ for }q\leq
2381,\,q=2393,2399,2401,2417,2437; \\
t_{2}(2,q) &<&4.5\sqrt{q}-1\,~\mbox{ for }q\leq 2473,\,q=2503,2531,2549.
\end{eqnarray*}
\end{theorem}

For $2633\leq q\leq 5107$ and for a few sporadic $q>5107$, the
values of $ \overline{t}_{2}(2,q)$ (up to June 2011) are
collected in \cite[Tabs\thinspace 3,4]{BDFMP-DM}. In this work
we obtained small arcs with new sizes for $q\in T_{6}.$ The new
arcs are obtained by computer search, based on the randomized
greedy algorithms. The current values of $\overline{t}
_{2}(2,q)<4.8\sqrt{q}$ for $2633\leq q\leq 5399$ are given in
Table 3. The data for $q\in T_{6}$ and other data (obtained in
this work) improving and extending results of
\cite[Tabs\thinspace 3,4]{BDFMP-DM} are written in Table~3 in
bold font. The data for $q=2659,2663,2683,2693,2753,2801$ with
$ \overline{t}_{2}(2,q)<4.5\sqrt{q}$ are written in italic
font.

\newpage

\begin{center}
\textbf{Table 3. }The smallest known sizes
$\overline{t}_{2}=\overline{t} _{2}(2,q)<4.8\sqrt{q}$ of
complete arcs in planes $PG(2,q),$ $2633\leq q\leq 5399$,
$A_{q}=\left\lfloor
a_{q}\sqrt{q}-\overline{t}_{2}(2,q)\right\rfloor $ , $B_{q}\geq
\overline{t}_{2}(2,q)/\sqrt{q}$\smallskip

$\renewcommand{\arraystretch}{0.9}
\begin{array}{@{}r@{\,\,\,\,}c@{\,\,\,\,}c@{\,\,\,}c|@{\,\,}c@{\,\,\,\,}c@{\,\,\,\,}c@{\,\,\,}c|@{\,\,}c@{\,\,\,\,}c@{\,\,\,\,}c@{\,\,\,}c|@{\,\,}c@{\,\,\,\,}c@{\,\,\,\,}c@{\,\,\,\,}c}
\hline
q & \overline{t}_{2} & A_{q} & B_{q} & q & \overline{t}_{2} & A_{q} & B_{q}
& q & \overline{t}_{2} & A_{q} & B_{q} & q & \overline{t}_{2} & A_{q} &
B_{q}^{\phantom{H^{L}}} \\ \hline
2633 & 231 & 25 & 4.51 & \mathbf{3307} & \mathbf{264} & \mathbf{23} &
\mathbf{4.60} & \mathbf{3947} & \mathbf{293} & \mathbf{21} & \mathbf{4.67} &
4663 & 323 & 18 & 4.74 \\
2647 & 232 & 25 & 4.51 & 3313 & 265 & 22 & 4.61 & 3967 & 294 & 20 & 4.67 &
\mathbf{4673} & \mathbf{323} & \mathbf{18} & \mathbf{4.73} \\
\mathbf{2657} & \mathbf{232} & \mathbf{25} & \mathbf{4.51} & 3319 & 265 & 23
& 4.60 & \mathbf{3989} & \mathbf{295} & \mathbf{20} & \mathbf{4.68} & 4679 &
324 & 18 & 4.74 \\
\mathbf{2659} & \emph{232} & \emph{0} & \emph{4.50} & 3323 & 265 & 23 & 4.60
& 4001 & 296 & 20 & 4.68 & 4691 & 324 & 18 & 4.74 \\
\mathbf{2663} & \emph{232} & \emph{0} & \emph{4.50} & \mathbf{3329} &
\mathbf{265} & \mathbf{23} & \mathbf{4.60} & 4003 & 296 & 20 & 4.68 & 4703 &
325 & 17 & 4.74 \\
2671 & 233 & 25 & 4.51 & \mathbf{3331} & \mathbf{265} & \mathbf{23} &
\mathbf{4.60} & \mathbf{4007} & \mathbf{296} & \mathbf{20} & \mathbf{4.68} &
4721 & 326 & 17 & 4.75 \\
\mathbf{2677} & \mathbf{233} & \mathbf{25} & \mathbf{4.51} & 3343 & 265 & 24
& 4.59 & 4013 & 296 & 20 & 4.68 & 4723 & 326 & 17 & 4.75 \\
\mathbf{2683} & \emph{233} & \emph{0} & \emph{4.50} & 3347 & 266 & 23 & 4.60
& 4019 & 296 & 20 & 4.67 & \mathbf{4729} & \mathbf{325} & \mathbf{18} &
\mathbf{4.73} \\
2687 & 234 & 25 & 4.52 & 3359 & 267 & 22 & 4.61 & 4021 & 296 & 21 & 4.67 &
4733 & 326 & 17 & 4.74 \\
2689 & 234 & 25 & 4.52 & 3361 & 267 & 22 & 4.61 & 4027 & 296 & 21 & 4.67 &
\mathbf{4751} & \mathbf{327} & \mathbf{17} & \mathbf{4.75} \\
\emph{2693} & \emph{233} & \emph{0} & \emph{4.49} & \mathbf{3371} & \mathbf{
267} & \mathbf{23} & \mathbf{4.60} & 4049 & 298 & 20 & 4.69 & 4759 & 327 & 17
& 4.75 \\
\mathbf{2699} & \mathbf{234} & \mathbf{25} & \mathbf{4.51} & \mathbf{3373} &
\mathbf{267} & \mathbf{23} & \mathbf{4.60} & \mathbf{4051} & \mathbf{298} &
\mathbf{20} & \mathbf{4.69} & 4783 & 328 & 17 & 4.75 \\
2707 & 235 & 25 & 4.52 & 3389 & 268 & 23 & 4.61 & 4057 & 298 & 20 & 4.68 &
4787 & 329 & 16 & 4.76 \\
2711 & 235 & 25 & 4.52 & \mathbf{3391} & \mathbf{267} & \mathbf{24} &
\mathbf{4.59} & 4073 & 299 & 20 & 4.69 & 4789 & 329 & 17 & 4.76 \\
2713 & 235 & 25 & 4.52 & \mathbf{3407} & \mathbf{269} & \mathbf{22} &
\mathbf{4.61} & \mathbf{4079} & \mathbf{299} & \mathbf{20} & \mathbf{4.69} &
\mathbf{4793} & \mathbf{329} & \mathbf{17} & \mathbf{4.76} \\
\mathbf{2719} & \mathbf{235} & \mathbf{25} & \mathbf{4.51} & 3413 & 269 & 23
& 4.61 & 4091 & 300 & 19 & 4.70 & \mathbf{4799} & \mathbf{329} & \mathbf{17}
& \mathbf{4.75} \\
2729 & 236 & 25 & 4.52 & 3433 & 270 & 22 & 4.61 & 4093 & 300 & 19 & 4.69 &
4801 & 329 & 17 & 4.75 \\
2731 & 236 & 25 & 4.52 & \mathbf{3449} & \mathbf{271} & \mathbf{22} &
\mathbf{4.62} & \mathbf{4096} & \mathbf{300} & \mathbf{20} & \mathbf{4.69} &
4813 & 330 & 16 & 4.76 \\
\mathbf{2741} & \mathbf{236} & \mathbf{25} & \mathbf{4.51} & 3457 & 271 & 22
& 4.61 & 4099 & 300 & 20 & 4.69 & 4817 & 329 & 18 & 4.75 \\
2749 & 237 & 25 & 4.53 & \mathbf{3461} & \mathbf{271} & \mathbf{23} &
\mathbf{4.61} & 4111 & 301 & 19 & 4.70 & 4831 & 329 & 18 & 4.74 \\
\emph{2753} & \emph{236} & \emph{0} & \emph{4.50} & 3463 & 272 & 22 & 4.63 &
\mathbf{4127} & \mathbf{301} & \mathbf{20} & \mathbf{4.69} & 4861 & 331 & 17
& 4.75 \\
2767 & 238 & 25 & 4.53 & 3467 & 272 & 22 & 4.62 & \mathbf{4129} & \mathbf{301
} & \mathbf{20} & \mathbf{4.69} & 4871 & 332 & 16 & 4.76 \\
2777 & 238 & 25 & 4.52 & 3469 & 272 & 22 & 4.62 & 4133 & 302 & 19 & 4.70 &
4877 & 332 & 17 & 4.76 \\
2789 & 239 & 25 & 4.53 & 3481 & 272 & 23 & 4.62 & 4139 & 302 & 19 & 4.70 &
4889 & 333 & 16 & 4.77 \\
2791 & 239 & 25 & 4.53 & 3491 & 273 & 22 & 4.63 & 4153 & 301 & 21 & 4.68 &
\mathbf{4903} & \mathbf{333} & \mathbf{17} & \mathbf{4.76} \\
\mathbf{2797} & \mathbf{239} & \mathbf{25} & \mathbf{4.52} & 3499 & 273 & 22
& 4.62 & 4157 & 303 & 19 & 4.70 & 4909 & 334 & 16 & 4.77 \\
\mathbf{2801} & \emph{238} & \emph{0} & \emph{4.50} & 3511 & 274 & 22 & 4.63
& 4159 & 302 & 20 & 4.69 & 4913 & 334 & 16 & 4.77 \\
2803 & 240 & 24 & 4.54 & 3517 & 274 & 22 & 4.63 & 4177 & 303 & 20 & 4.69 &
4919 & 334 & 16 & 4.77 \\
2809 & 240 & 25 & 4.53 & \mathbf{3527} & \mathbf{274} & \mathbf{22} &
\mathbf{4.62} & \mathbf{4201} & \mathbf{304} & \mathbf{20} & \mathbf{4.70} &
\mathbf{4931} & \mathbf{334} & \mathbf{17} & \mathbf{4.76} \\
\mathbf{2819} & \mathbf{240} & \mathbf{25} & \mathbf{4.53} & 3529 & 275 & 22
& 4.63 & 4211 & 305 & 19 & 4.71 & 4933 & 335 & 16 & 4.77 \\
\mathbf{2833} & \mathbf{240} & \mathbf{26} & \mathbf{4.51} & 3533 & 275 & 22
& 4.63 & 4217 & 305 & 19 & 4.70 & 4937 & 335 & 16 & 4.77 \\
\mathbf{2837} & \mathbf{241} & \mathbf{25} & \mathbf{4.53} & 3539 & 275 & 22
& 4.63 & 4219 & 305 & 19 & 4.70 & 4943 & 334 & 17 & 4.76 \\
2843 & 242 & 24 & 4.54 & \mathbf{3541} & \mathbf{275} & \mathbf{22} &
\mathbf{4.63} & 4229 & 305 & 20 & 4.70 & 4951 & 335 & 16 & 4.77 \\
\mathbf{2851} & \mathbf{242} & \mathbf{24} & \mathbf{4.54} & \mathbf{3547} &
\mathbf{275} & \mathbf{22} & \mathbf{4.62} & 4231 & 306 & 19 & 4.71 & 4957 &
335 & 17 & 4.76 \\
\mathbf{2857} & \mathbf{242} & \mathbf{25} & \mathbf{4.53} & \mathbf{3557} &
\mathbf{276} & \mathbf{22} & \mathbf{4.63} & 4241 & 306 & 19 & 4.70 & 4967 &
336 & 16 & 4.77 \\
2861 & 243 & 24 & 4.55 & 3559 & 276 & 22 & 4.63 & 4243 & 306 & 19 & 4.70 &
4969 & 336 & 16 & 4.77 \\
\mathbf{2879} & \mathbf{243} & \mathbf{25} & \mathbf{4.53} & 3571 & 277 & 21
& 4.64 & 4253 & 306 & 20 & 4.70 & \mathbf{4973} & \mathbf{336} & \mathbf{16}
& \mathbf{4.77} \\
2887 & 243 & 25 & 4.53 & \mathbf{3581} & \mathbf{277} & \mathbf{22} &
\mathbf{4.63} & 4259 & 307 & 19 & 4.71 & 4987 & 336 & 17 & 4.76 \\
\mathbf{2897} & \mathbf{244} & \mathbf{25} & \mathbf{4.54} & 3583 & 277 & 22
& 4.63 & 4261 & 307 & 19 & 4.71 & 4993 & 337 & 16 & 4.77 \\
2903 & 245 & 24 & 4.55 & 3593 & 278 & 21 & 4.64 & 4271 & 307 & 19 & 4.70 &
\mathbf{4999} & \mathbf{337} & \mathbf{16} & \mathbf{4.77} \\
2909 & 245 & 24 & 4.55 & 3607 & 278 & 22 & 4.63 & 4273 & 306 & 20 & 4.69 &
5003 & 337 & 16 & 4.77 \\
\mathbf{2917} & \mathbf{245} & \mathbf{25} & \mathbf{4.54} & \mathbf{3613} &
\mathbf{278} & \mathbf{22} & \mathbf{4.63} & 4283 & 308 & 19 & 4.71 & 5009 &
337 & 16 & 4.77 \\ \hline
\end{array}
$\newpage

\textbf{Table 3} (continue)\smallskip

$\renewcommand{\arraystretch}{0.9}
\begin{array}{@{}r@{\,\,\,\,}c@{\,\,\,\,}c@{\,\,\,}c|@{\,\,}c@{\,\,\,\,}c@{\,\,\,\,}c@{\,\,\,}c|@{\,\,}c@{\,\,\,\,}c@{\,\,\,\,}c@{\,\,\,}c|@{\,\,}c@{\,\,\,\,}c@{\,\,\,\,}c@{\,\,\,\,}c}
\hline
q & \overline{t}_{2} & A_{q} & B_{q} & q & \overline{t}_{2} & A_{q} & B_{q}
& q & \overline{t}_{2} & A_{q} & B_{q} & q & \overline{t}_{2} & A_{q} &
B_{q}^{\phantom{H^{L}}} \\ \hline
2927 & 245 & 25 & 4.53 & 3617 & 278 & 22 & 4.63 & 4289 & 308 & 19 & 4.71 &
5011 & 337 & 16 & 4.77 \\
2939 & 246 & 25 & 4.54 & 3623 & 279 & 21 & 4.64 & 4297 & 308 & 19 & 4.70 &
5021 & 338 & 16 & 4.78 \\
\mathbf{2953} & \mathbf{246} & \mathbf{25} & \mathbf{4.53} & \mathbf{3631} &
\mathbf{279} & \mathbf{22} & \mathbf{4.64} & 4327 & 310 & 18 & 4.72 &
\mathbf{5023} & \mathbf{338} & \mathbf{16} & \mathbf{4.77} \\
\mathbf{2957} & \mathbf{247} & \mathbf{24} & \mathbf{4.55} & 3637 & 280 & 21
& 4.65 & \mathbf{4337} & \mathbf{310} & \mathbf{19} & \mathbf{4.71} & 5039 &
339 & 15 & 4.78 \\
2963 & 248 & 24 & 4.56 & 3643 & 278 & 23 & 4.61 & \mathbf{4339} & \mathbf{310
} & \mathbf{19} & \mathbf{4.71} & 5041 & 339 & 16 & 4.78 \\
2969 & 248 & 24 & 4.56 & 3659 & 281 & 21 & 4.65 & 4349 & 311 & 18 & 4.72 &
\mathbf{5051} & \mathbf{339} & \mathbf{16} & \mathbf{4.77} \\
\mathbf{2971} & \mathbf{247} & \mathbf{25} & \mathbf{4.54} & \mathbf{3671} &
\mathbf{281} & \mathbf{21} & \mathbf{4.64} & 4357 & 311 & 19 & 4.72 & 5059 &
339 & 16 & 4.77 \\
\mathbf{2999} & \mathbf{249} & \mathbf{24} & \mathbf{4.55} & \mathbf{3673} &
\mathbf{280} & \mathbf{23} & \mathbf{4.63} & 4363 & 310 & 20 & 4.70 &
\mathbf{5077} & \mathbf{340} & \mathbf{16} & \mathbf{4.78} \\
3001 & 250 & 23 & 4.57 & \mathbf{3677} & \mathbf{281} & \mathbf{22} &
\mathbf{4.64} & 4373 & 312 & 18 & 4.72 & \mathbf{5081} & \mathbf{340} &
\mathbf{16} & \mathbf{4.77} \\
\mathbf{3011} & \mathbf{250} & \mathbf{24} & \mathbf{4.56} & \mathbf{3691} &
\mathbf{282} & \mathbf{21} & \mathbf{4.65} & \mathbf{4391} & \mathbf{312} &
\mathbf{19} & \mathbf{4.71} & 5087 & 341 & 15 & 4.79 \\
\mathbf{3019} & \mathbf{250} & \mathbf{24} & \mathbf{4.55} & \mathbf{3697} &
\mathbf{282} & \mathbf{22} & \mathbf{4.64} & 4397 & 313 & 18 & 4.73 &
\mathbf{5099} & \mathbf{341} & \mathbf{16} & \mathbf{4.78} \\
3023 & 251 & 23 & 4.57 & \mathbf{3701} & \mathbf{282} & \mathbf{22} &
\mathbf{4.64} & \mathbf{4409} & \mathbf{313} & \mathbf{19} & \mathbf{4.72} &
\mathbf{5101} & \mathbf{341} & \mathbf{16} & \mathbf{4.78} \\
\mathbf{3037} & \mathbf{251} & \mathbf{24} & \mathbf{4.56} & 3709 & 283 & 21
& 4.65 & 4421 & 314 & 18 & 4.73 & 5107 & 341 & 16 & 4.78 \\
\mathbf{3041} & \mathbf{251} & \mathbf{24} & \mathbf{4.56} & \mathbf{3719} &
\mathbf{283} & \mathbf{21} & \mathbf{4.65} & 4423 & 312 & 20 & 4.70 &
\mathbf{5113} & \mathbf{341} & \mathbf{16} & \mathbf{4.77} \\
3049 & 252 & 24 & 4.57 & \mathbf{3721} & \mathbf{283} & \mathbf{22} &
\mathbf{4.64} & 4441 & 315 & 18 & 4.73 & 5119 & 341 & 16 & 4.77 \\
\mathbf{3061} & \mathbf{252} & \mathbf{24} & \mathbf{4.56} & 3727 & 284 & 21
& 4.66 & 4447 & 314 & 19 & 4.71 & 5147 & 343 & 15 & 4.79 \\
3067 & 253 & 23 & 4.57 & 3733 & 284 & 21 & 4.65 & \mathbf{4451} & \mathbf{315
} & \mathbf{18} & \mathbf{4.73} & \mathbf{5153} & \mathbf{342} & \mathbf{16}
& \mathbf{4.77} \\
3079 & 253 & 24 & 4.56 & 3739 & 283 & 22 & 4.63 & 4457 & 315 & 18 & 4.72 &
\mathbf{5167} & \mathbf{344} & \mathbf{15} & \mathbf{4.79} \\
3083 & 253 & 24 & 4.56 & \mathbf{3761} & \mathbf{284} & \mathbf{22} &
\mathbf{4.64} & 4463 & 315 & 19 & 4.72 & \mathbf{5171} & \mathbf{344} &
\mathbf{15} & \mathbf{4.79} \\
3089 & 254 & 23 & 4.58 & \mathbf{3767} & \mathbf{285} & \mathbf{21} &
\mathbf{4.65} & 4481 & 315 & 19 & 4.71 & \mathbf{5179} & \mathbf{344} &
\mathbf{15} & \mathbf{4.79} \\
3109 & 255 & 23 & 4.58 & 3769 & 286 & 20 & 4.66 & \mathbf{4483} & \mathbf{316
} & \mathbf{18} & \mathbf{4.72} & \mathbf{5189} & \mathbf{344} & \mathbf{16}
& \mathbf{4.78} \\
3119 & 255 & 24 & 4.57 & 3779 & 286 & 21 & 4.66 & 4489 & 316 & 19 & 4.72 &
\mathbf{5197} & \mathbf{345} & \mathbf{15} & \mathbf{4.79} \\
3121 & 255 & 24 & 4.57 & 3793 & 287 & 20 & 4.67 & 4493 & 317 & 18 & 4.73 &
\mathbf{5209} & \mathbf{345} & \mathbf{15} & \mathbf{4.79} \\
3125 & 256 & 23 & 4.58 & 3797 & 287 & 21 & 4.66 & \mathbf{4507} & \mathbf{317
} & \mathbf{18} & \mathbf{4.73} & \mathbf{5227} & \mathbf{346} & \mathbf{15}
& \mathbf{4.79} \\
\mathbf{3137} & \mathbf{256} & \mathbf{24} & \mathbf{4.58} & 3803 & 287 & 21
& 4.66 & 4513 & 318 & 17 & 4.74 & \mathbf{5231} & \mathbf{346} & \mathbf{15}
& \mathbf{4.79} \\
3163 & 257 & 24 & 4.57 & 3821 & 288 & 21 & 4.66 & 4517 & 317 & 19 & 4.72 &
\mathbf{5233} & \mathbf{346} & \mathbf{15} & \mathbf{4.79} \\
3167 & 258 & 23 & 4.59 & \mathbf{3823} & \mathbf{288} & \mathbf{21} &
\mathbf{4.66} & 4519 & 318 & 18 & 4.74 & 5237 & 347 & 14 & 4.80 \\
3169 & 258 & 23 & 4.59 & \mathbf{3833} & \mathbf{288} & \mathbf{21} &
\mathbf{4.66} & 4523 & 318 & 18 & 4.73 & \mathbf{5261} & \mathbf{347} &
\mathbf{15} & \mathbf{4.79} \\
\mathbf{3181} & \mathbf{258} & \mathbf{24} & \mathbf{4.58} & \mathbf{3847} &
\mathbf{288} & \mathbf{22} & \mathbf{4.65} & 4547 & 319 & 18 & 4.74 &
\mathbf{5273} & \mathbf{348} & \mathbf{15} & \mathbf{4.80} \\
3187 & 258 & 24 & 4.58 & \mathbf{3851} & \mathbf{288} & \mathbf{22} &
\mathbf{4.65} & 4549 & 319 & 18 & 4.73 & \mathbf{5279} & \mathbf{348} &
\mathbf{15} & \mathbf{4.79} \\
3191 & 259 & 23 & 4.59 & 3853 & 289 & 21 & 4.66 & 4561 & 319 & 18 & 4.73 &
\mathbf{5281} & \mathbf{348} & \mathbf{15} & \mathbf{4.79} \\
3203 & 259 & 23 & 4.58 & 3863 & 290 & 20 & 4.67 & 4567 & 320 & 17 & 4.74 &
\mathbf{5297} & \mathbf{349} & \mathbf{14} & \mathbf{4.80} \\
3209 & 260 & 23 & 4.59 & \mathbf{3877} & \mathbf{290} & \mathbf{21} &
\mathbf{4.66} & 4583 & 320 & 18 & 4.73 & 5303 & 349 & 15 & 4.80 \\
\mathbf{3217} & \mathbf{260} & \mathbf{23} & \mathbf{4.59} & 3881 & 291 & 20
& 4.68 & 4591 & 321 & 17 & 4.74 & \mathbf{5309} & \mathbf{349} & \mathbf{15}
& \mathbf{4.79} \\
\mathbf{3221} & \mathbf{260} & \mathbf{23} & \mathbf{4.59} & 3889 & 291 & 20
& 4.67 & 4597 & 321 & 18 & 4.74 & \mathbf{5323} & \mathbf{349} & \mathbf{15}
& \mathbf{4.79} \\
3229 & 260 & 24 & 4.58 & 3907 & 292 & 20 & 4.68 & \mathbf{4603} & \mathbf{321
} & \mathbf{18} & \mathbf{4.74} & \mathbf{5329} & \mathbf{350} & \mathbf{15}
& \mathbf{4.80} \\
3251 & 262 & 23 & 4.60 & 3911 & 292 & 20 & 4.67 & \mathbf{4621} & \mathbf{321
} & \mathbf{18} & \mathbf{4.73} & \mathbf{5333} & \mathbf{349} & \mathbf{16}
& \mathbf{4.78} \\
3253 & 261 & 24 & 4.58 & \mathbf{3917} & \mathbf{291} & \mathbf{21} &
\mathbf{4.65} & 4637 & 322 & 18 & 4.73 & \mathbf{5347} & \mathbf{350} &
\mathbf{15} & \mathbf{4.79} \\
3257 & 262 & 23 & 4.60 & 3919 & 292 & 21 & 4.67 & 4639 & 323 & 17 & 4.75 &
\mathbf{5351} & \mathbf{350} & \mathbf{15} & \mathbf{4.79} \\
\mathbf{3259} & \mathbf{262} & \mathbf{23} & \mathbf{4.59} & \mathbf{3923} &
\mathbf{292} & \mathbf{21} & \mathbf{4.67} & 4643 & 323 & 17 & 4.75 &
\mathbf{5381} & \mathbf{352} & \mathbf{14} & \mathbf{4.80} \\
3271 & 263 & 22 & 4.60 & 3929 & 293 & 20 & 4.68 & 4649 & 323 & 17 & 4.74 &
\mathbf{5387} & \mathbf{352} & \mathbf{14} & \mathbf{4.80} \\
3299 & 264 & 23 & 4.60 & 3931 & 293 & 20 & 4.68 & 4651 & 323 & 17 & 4.74 &
\mathbf{5393} & \mathbf{352} & \mathbf{15} & \mathbf{4.80} \\
3301 & 264 & 23 & 4.60 & \mathbf{3943} & \mathbf{293} & \mathbf{20} &
\mathbf{4.67} & 4657 & 323 & 18 & 4.74 & \mathbf{5399} & \mathbf{352} &
\mathbf{15} & \mathbf{4.80} \\ \hline
\end{array}
$\newpage
\end{center}

\newpage

The current values of $\overline{t}_{2}(2,q)$ for $5407\leq
q\leq 8353$ and $ 8363\leq q\leq 9109$ are given in Tables 4
and 5, respectively. All results in these tables are new and
they have been obtained in this work by computer search, based
on the randomized greedy algorithms. Data for $
q=5413,5417,5419,5441,5443,5471,5483,5501,5521$ with
$\overline{t} _{2}(2,q)<4.8\sqrt{q}$ are written in Table~4 in
italic font.

From Tables 3 - 5 we obtain Theorem \ref{Th3_5sqroot(q)} improving and
extending the results of \cite[Th.\thinspace 3.3]{BDFMP-DM}.

\begin{theorem}
\label{Th3_5sqroot(q)} In $PG(2,q),$ the following holds.
\begin{eqnarray}
t_{2}(2,q) &<&5\sqrt{q}~~\,\mbox{ for }q\leq 9067. \\
t_{2}(2,q) &<&4.6\sqrt{q}~~\mbox{ for }q\leq
3307,\,q=3319,3323,3329,3331,3343,3347,3371,3373,3391;  \notag \\
t_{2}(2,q) &<&4.7\sqrt{q}~~\mbox{ for }q\leq
4201,\,q=4217,4219,4229,4241,4243,4253,4271,4273,4297,  \notag \\
&&\phantom{4.7\sqrt{q}~~\mbox{ for }q\leq 4057,\,q=}~4363,4423.  \notag \\
t_{2}(2,q) &<&4.8\sqrt{q}~~\mbox{ for }q\leq
5399,\,q=5413,5417,5419,5441,5443,5471,5483,5501,5521;  \notag \\
t_{2}(2,q) &<&4.9\sqrt{q}~~\mbox{ for }q\leq
6907,\,q= 6947,6949,6961,6971,6983,6997,7001,7039,7187, \notag \\
&&\phantom{4.7\sqrt{q}~~\mbox{ for }q\leq 4057,\,q=}
~7193,7307,7451.  \notag
\end{eqnarray}
Also,
\begin{eqnarray*}
t_{2}(2,q) &<&5\sqrt{q}-22\mbox{ for }q\leq
3559,\,q=3581,3583,3607,3613,3617,3631,3643,3673,3677, \\
&&\phantom{5\sqrt{q}-22\,~ \mbox{ for } q\leq 1369, q=}
3697,3701,3721,3739,3761,3847,3851; \\
t_{2}(2,q) &<&5\sqrt{q}-21\mbox{ for }q\leq
3767,\,q=3779,3797,3803,3821,3823,3833,3847,3851,3853, \\
&&\phantom{5\sqrt{q}-22\,~ \mbox{ for } q\leq 1369, q=}
3877,3917,3919,3923,3947,4021,4027,4153; \\
t_{2}(2,q) &<&5\sqrt{q}-20\mbox{ for }q\leq
4079,\,q=4096,4099,4127,4129,4153,4159,4177,4201,4229, \\
&&\phantom{5\sqrt{q}-20\mbox{ for }q\leq 3911,\,q=\,\,}4253,4273,4363,4423;
\\
t_{2}(2,q) &<&5\sqrt{q}-19\mbox{ for }q\leq
4297,\,q=4337,4339,4357,4363,4391,4409,4423,4447,4463, \\
&&\phantom{5\sqrt{q}-20\mbox{ for }q\leq 3911,\,q=\,\,}4481,4489,4517; \\
t_{2}(2,q) &<&5\sqrt{q}-16\mbox{ for }q\leq
5023,\,q=5041,5051,5059,5077,5081,5099,5101,5107,5113, \\
&&\phantom{5\sqrt{q}-20\mbox{ for }q\leq 3911,\,q=\,\,}5119,5153,5189,5333;
\\
t_{2}(2,q) &<&5\sqrt{q}-14\mbox{ for }q\leq
5501,\,q=5507,5519,5521,5527,5557,5569,5573,5581,5591, \\
&&\phantom{5\sqrt{q}-20\mbox{ for }q\leq 3911,\,q=\,\,}
5689,5693,5711,5717,5749,5783,5813; \\
t_{2}(2,q) &<&5\sqrt{q}-12\mbox{ for }q\leq
5881,\,q=5903,5923,5927,5939,5953,5987,6007,6029,6053, \\
&&\phantom{5\sqrt{q}-16\mbox{ for }q\leq 4969,\,q=\,\,}
6073,6089,6143,6151,6163.
\end{eqnarray*}
\newpage
\end{theorem}

\begin{center}
\textbf{Table 4. }The smallest known sizes
$\overline{t}_{2}=\overline{t} _{2}(2,q)$ of complete arcs in
planes $PG(2,q),$ $5407\leq q\leq 8353$, $ A_{q}=\left\lfloor
5\sqrt{q}-\overline{t}_{2}(2,q)\right\rfloor $, $ B_{q}\geq
\overline{t}_{2}(2,q)/\sqrt{q}$\smallskip

$\renewcommand{\arraystretch}{0.87}
\begin{array}{@{}r@{\,\,\,\,}c@{\,\,\,\,}c@{\,\,\,}c|@{\,\,}c@{\,\,\,\,}c@{\,\,\,\,}c@{\,\,\,}c|@{\,\,}c@{\,\,\,\,}c@{\,\,\,\,}c@{\,\,\,}c|@{\,\,}c@{\,\,\,\,}c@{\,\,\,\,}c@{\,\,\,\,}c}
\hline
q & \overline{t}_{2} & A_{q} & B_{q} & q & \overline{t}_{2} & A_{q} & B_{q}
& q & \overline{t}_{2} & A_{q} & B_{q} & q & \overline{t}_{2} & A_{q} &
B_{q}^{\phantom{H^{L}}} \\ \hline
5407 & 353 & 14 & 4.81 & 6121 & 380 & 11 & 4.86 & 6841 & 404 & 9 & 4.89 &
7589 & 430 & 5 & 4.94 \\
\emph{5413} & \emph{353} & \emph{14} & \emph{4.80} & 6131 & 380 & 11 & 4.86
& 6857 & 405 & 9 & 4.90 & 7591 & 430 & 5 & 4.94 \\
\emph{5417} & \emph{353} & \emph{15} & \emph{4.80} & 6133 & 380 & 11 & 4.86
& 6859 & 405 & 9 & 4.90 & 7603 & 430 & 5 & 4.94 \\
\emph{5419} & \emph{353} & \emph{15} & \emph{4.80} & 6143 & 379 & 12 & 4.84
& 6863 & 405 & 9 & 4.89 & 7607 & 430 & 6 & 4.94 \\
5431 & 354 & 14 & 4.81 & 6151 & 380 & 12 & 4.85 & 6869 & 405 & 9 & 4.89 &
7621 & 431 & 5 & 4.94 \\
5437 & 354 & 14 & 4.81 & 6163 & 380 & 12 & 4.85 & 6871 & 405 & 9 & 4.89 &
7639 & 432 & 5 & 4.95 \\
\emph{5441} & \emph{354} & \emph{14} & \emph{4.80} & 6173 & 382 & 10 & 4.87
& 6883 & 406 & 8 & 4.90 & 7643 & 431 & 6 & 4.93 \\
\emph{5443} & \emph{354} & \emph{14} & \emph{4.80} & 6197 & 382 & 11 & 4.86
& 6889 & 406 & 9 & 4.90 & 7649 & 432 & 5 & 4.94 \\
5449 & 355 & 14 & 4.81 & 6199 & 383 & 10 & 4.87 & 6899 & 406 & 9 & 4.89 &
7669 & 431 & 6 & 4.93 \\
\emph{5471} & \emph{355} & \emph{14} & \emph{4.80} & 6203 & 383 & 10 & 4.87
& 6907 & 407 & 8 & 4.90 & 7673 & 432 & 5 & 4.94 \\
5477 & 356 & 14 & 4.82 & 6211 & 383 & 11 & 4.86 & 6911 & 408 & 7 & 4.91 &
7681 & 433 & 5 & 4.95 \\
5479 & 356 & 14 & 4.81 & 6217 & 383 & 11 & 4.86 & 6917 & 408 & 7 & 4.91 &
7687 & 433 & 5 & 4.94 \\
\emph{5483} & \emph{355} & \emph{15} & \emph{4.80} & 6221 & 383 & 11 & 4.86
& 6947 & 408 & 8 & 4.90 & 7691 & 433 & 5 & 4.94 \\
\emph{5501} & \emph{356} & \emph{14} & \emph{4.80} & 6229 & 383 & 11 & 4.86
& 6949 & 408 & 8 & 4.90 & 7699 & 433 & 5 & 4.94 \\
5503 & 357 & 13 & 4.82 & 6241 & 384 & 11 & 4.87 & 6959 & 409 & 8 & 4.91 &
7703 & 434 & 4 & 4.95 \\
5507 & 357 & 14 & 4.82 & 6247 & 384 & 11 & 4.86 & 6961 & 408 & 9 & 4.90 &
7717 & 434 & 5 & 4.95 \\
5519 & 357 & 14 & 4.81 & 6257 & 384 & 11 & 4.86 & 6967 & 409 & 8 & 4.91 &
7723 & 435 & 4 & 4.95 \\
5521 & 356 & 15 & 4.80 & 6263 & 385 & 10 & 4.87 & 6971 & 409 & 8 & 4.90 &
7727 & 434 & 5 & 4.94 \\
5527 & 357 & 14 & 4.81 & 6269 & 384 & 11 & 4.85 & 6977 & 410 & 7 & 4.91 &
7741 & 435 & 4 & 4.95 \\
5531 & 358 & 13 & 4.82 & 6271 & 384 & 11 & 4.85 & 6983 & 408 & 9 & 4.89 &
7753 & 436 & 4 & 4.96 \\
5557 & 358 & 14 & 4.81 & 6277 & 385 & 11 & 4.86 & 6991 & 410 & 8 & 4.91 &
7757 & 436 & 4 & 4.96 \\
5563 & 359 & 13 & 4.82 & 6287 & 385 & 11 & 4.86 & 6997 & 409 & 9 & 4.89 &
7759 & 436 & 4 & 4.95 \\
5569 & 359 & 14 & 4.82 & 6299 & 385 & 11 & 4.86 & 7001 & 409 & 9 & 4.89 &
7789 & 437 & 4 & 4.96 \\
5573 & 359 & 14 & 4.81 & 6301 & 386 & 10 & 4.87 & 7013 & 411 & 7 & 4.91 &
7793 & 437 & 4 & 4.96 \\
5581 & 359 & 14 & 4.81 & 6311 & 386 & 11 & 4.86 & 7019 & 411 & 7 & 4.91 &
7817 & 437 & 5 & 4.95 \\
5591 & 359 & 14 & 4.81 & 6317 & 387 & 10 & 4.87 & 7027 & 412 & 7 & 4.92 &
7823 & 438 & 4 & 4.96 \\
5623 & 361 & 13 & 4.82 & 6323 & 387 & 10 & 4.87 & 7039 & 410 & 9 & 4.89 &
7829 & 437 & 5 & 4.94 \\
5639 & 362 & 13 & 4.83 & 6329 & 387 & 10 & 4.87 & 7043 & 412 & 7 & 4.91 &
7841 & 438 & 4 & 4.95 \\
5641 & 362 & 13 & 4.82 & 6337 & 388 & 10 & 4.88 & 7057 & 413 & 7 & 4.92 &
7853 & 438 & 5 & 4.95 \\
5647 & 362 & 13 & 4.82 & 6343 & 388 & 10 & 4.88 & 7069 & 412 & 8 & 4.91 &
7867 & 440 & 3 & 4.97 \\
5651 & 362 & 13 & 4.82 & 6353 & 388 & 10 & 4.87 & 7079 & 413 & 7 & 4.91 &
7873 & 440 & 3 & 4.96 \\
5653 & 362 & 13 & 4.82 & 6359 & 387 & 11 & 4.86 & 7103 & 414 & 7 & 4.92 &
7877 & 438 & 5 & 4.94 \\
5657 & 363 & 13 & 4.83 & 6361 & 388 & 10 & 4.87 & 7109 & 414 & 7 & 4.92 &
7879 & 440 & 3 & 4.96 \\
5659 & 363 & 13 & 4.83 & 6367 & 389 & 9 & 4.88 & 7121 & 414 & 7 & 4.91 & 7883
& 438 & 5 & 4.94 \\
5669 & 363 & 13 & 4.83 & 6373 & 389 & 10 & 4.88 & 7127 & 414 & 8 & 4.91 &
7901 & 441 & 3 & 4.97 \\
5683 & 363 & 13 & 4.82 & 6379 & 389 & 10 & 4.88 & 7129 & 415 & 7 & 4.92 &
7907 & 440 & 4 & 4.95 \\
5689 & 363 & 14 & 4.82 & 6389 & 388 & 11 & 4.86 & 7151 & 416 & 6 & 4.92 &
7919 & 441 & 3 & 4.96 \\
5693 & 363 & 14 & 4.82 & 6397 & 389 & 10 & 4.87 & 7159 & 415 & 8 & 4.91 &
7921 & 440 & 5 & 4.95 \\
5701 & 364 & 13 & 4.83 & 6421 & 390 & 10 & 4.87 & 7177 & 416 & 7 & 4.92 &
7927 & 440 & 5 & 4.95 \\
5711 & 363 & 14 & 4.81 & 6427 & 391 & 9 & 4.88 & 7187 & 415 & 8 & 4.90 & 7933
& 442 & 3 & 4.97 \\
5717 & 364 & 14 & 4.82 & 6449 & 391 & 10 & 4.87 & 7193 & 415 & 9 & 4.90 &
7937 & 442 & 3 & 4.97 \\
5737 & 365 & 13 & 4.82 & 6451 & 392 & 9 & 4.89 & 7207 & 417 & 7 & 4.92 & 7949
& 442 & 3 & 4.96 \\ \hline
\end{array}
$ \newpage

\textbf{Table 4 }(continue)\smallskip

$\renewcommand{\arraystretch}{0.89}
\begin{array}{@{}r@{\,\,\,\,}c@{\,\,\,\,}c@{\,\,\,}c|@{\,\,}c@{\,\,\,\,}c@{\,\,\,\,}c@{\,\,\,}c|@{\,\,}c@{\,\,\,\,}c@{\,\,\,\,}c@{\,\,\,}c|@{\,\,}c@{\,\,\,\,}c@{\,\,\,\,}c@{\,\,\,\,}c}
\hline
q & \overline{t}_{2} & A_{q} & B_{q} & q & \overline{t}_{2} & A_{q} & B_{q}
& q & \overline{t}_{2} & A_{q} & B_{q} & q & \overline{t}_{2} & A_{q} &
B_{q}^{\phantom{H^{L}}} \\ \hline
5741 & 366 & 12 & 4.84 & 6469 & 391 & 11 & 4.87 & 7211 & 418 & 6 & 4.93 &
7951 & 441 & 4 & 4.95 \\
5743 & 365 & 13 & 4.82 & 6473 & 392 & 10 & 4.88 & 7213 & 417 & 7 & 4.91 &
7963 & 442 & 4 & 4.96 \\
5749 & 365 & 14 & 4.82 & 6481 & 392 & 10 & 4.87 & 7219 & 417 & 7 & 4.91 &
7993 & 443 & 4 & 4.96 \\
5779 & 367 & 13 & 4.83 & 6491 & 393 & 9 & 4.88 & 7229 & 418 & 7 & 4.92 & 8009
& 444 & 3 & 4.97 \\
5783 & 366 & 14 & 4.82 & 6521 & 392 & 11 & 4.86 & 7237 & 418 & 7 & 4.92 &
8011 & 444 & 3 & 4.97 \\
5791 & 367 & 13 & 4.83 & 6529 & 393 & 11 & 4.87 & 7243 & 419 & 6 & 4.93 &
8017 & 443 & 4 & 4.95 \\
5801 & 367 & 13 & 4.82 & 6547 & 394 & 10 & 4.87 & 7247 & 419 & 6 & 4.93 &
8039 & 444 & 4 & 4.96 \\
5807 & 368 & 13 & 4.83 & 6551 & 395 & 9 & 4.89 & 7253 & 418 & 7 & 4.91 & 8053
& 445 & 3 & 4.96 \\
5813 & 366 & 15 & 4.81 & 6553 & 395 & 9 & 4.88 & 7283 & 420 & 6 & 4.93 & 8059
& 446 & 2 & 4.97 \\
5821 & 369 & 12 & 4.84 & 6561 & 395 & 10 & 4.88 & 7297 & 421 & 6 & 4.93 & 8069
& 446 & 3 & 4.97 \\
5827 & 369 & 12 & 4.84 & 6563 & 395 & 10 & 4.88 & 7307 & 418 & 9 & 4.89 &
8081 & 445 & 4 & 4.96 \\
5839 & 369 & 13 & 4.83 & 6569 & 395 & 10 & 4.88 & 7309 & 420 & 7 & 4.92 &
8087 & 446 & 3 & 4.96 \\
5843 & 370 & 12 & 4.85 & 6571 & 396 & 9 & 4.89 & 7321 & 421 & 6 & 4.93 & 8089
& 447 & 2 & 4.98 \\
5849 & 369 & 13 & 4.83 & 6577 & 396 & 9 & 4.89 & 7331 & 422 & 6 & 4.93 & 8093
& 447 & 2 & 4.97 \\
5851 & 370 & 12 & 4.84 & 6581 & 396 & 9 & 4.89 & 7333 & 421 & 7 & 4.92 & 8101
& 447 & 3 & 4.97 \\
5857 & 370 & 12 & 4.84 & 6599 & 396 & 10 & 4.88 & 7349 & 422 & 6 & 4.93 &
8111 & 448 & 2 & 4.98 \\
5861 & 370 & 12 & 4.84 & 6607 & 397 & 9 & 4.89 & 7351 & 422 & 6 & 4.93 & 8117
& 448 & 2 & 4.98 \\
5867 & 370 & 12 & 4.84 & 6619 & 398 & 8 & 4.90 & 7369 & 422 & 7 & 4.92 & 8123
& 448 & 2 & 4.98 \\
5869 & 370 & 13 & 4.83 & 6637 & 398 & 9 & 4.89 & 7393 & 423 & 6 & 4.92 & 8147
& 448 & 3 & 4.97 \\
5879 & 371 & 12 & 4.84 & 6653 & 398 & 9 & 4.88 & 7411 & 425 & 5 & 4.94 & 8161
& 448 & 3 & 4.96 \\
5881 & 371 & 12 & 4.84 & 6659 & 398 & 10 & 4.88 & 7417 & 423 & 7 & 4.92 &
8167 & 448 & 3 & 4.96 \\
5897 & 372 & 11 & 4.85 & 6661 & 398 & 10 & 4.88 & 7433 & 425 & 6 & 4.93 &
8171 & 448 & 3 & 4.96 \\
5903 & 372 & 12 & 4.85 & 6673 & 399 & 9 & 4.89 & 7451 & 422 & 9 & 4.89 & 8179
& 449 & 3 & 4.97 \\
5923 & 372 & 12 & 4.84 & 6679 & 399 & 9 & 4.89 & 7457 & 424 & 7 & 4.92 & 8191
& 449 & 3 & 4.97 \\
5927 & 372 & 12 & 4.84 & 6689 & 399 & 9 & 4.88 & 7459 & 425 & 6 & 4.93 & 8192
& 449 & 3 & 4.97 \\
5939 & 373 & 12 & 4.85 & 6691 & 399 & 9 & 4.88 & 7477 & 426 & 6 & 4.93 & 8209
& 450 & 3 & 4.97 \\
5953 & 372 & 13 & 4.83 & 6701 & 400 & 9 & 4.89 & 7481 & 426 & 6 & 4.93 & 8219
& 451 & 2 & 4.98 \\
5981 & 375 & 11 & 4.85 & 6703 & 400 & 9 & 4.89 & 7487 & 426 & 6 & 4.93 & 8221
& 451 & 2 & 4.98 \\
5987 & 374 & 12 & 4.84 & 6709 & 400 & 9 & 4.89 & 7489 & 426 & 6 & 4.93 & 8231
& 451 & 2 & 4.98 \\
6007 & 375 & 12 & 4.84 & 6719 & 400 & 9 & 4.88 & 7499 & 427 & 5 & 4.94 & 8233
& 451 & 2 & 4.98 \\
6011 & 376 & 11 & 4.85 & 6733 & 401 & 9 & 4.89 & 7507 & 427 & 6 & 4.93 & 8237
& 452 & 1 & 4.99 \\
6029 & 375 & 13 & 4.83 & 6737 & 401 & 9 & 4.89 & 7517 & 427 & 6 & 4.93 & 8243
& 451 & 2 & 4.97 \\
6037 & 377 & 11 & 4.86 & 6761 & 402 & 9 & 4.89 & 7523 & 428 & 5 & 4.94 & 8263
& 451 & 3 & 4.97 \\
6043 & 377 & 11 & 4.85 & 6763 & 402 & 9 & 4.89 & 7529 & 428 & 5 & 4.94 & 8269
& 453 & 1 & 4.99 \\
6047 & 377 & 11 & 4.85 & 6779 & 403 & 8 & 4.90 & 7537 & 428 & 6 & 4.93 & 8273
& 452 & 2 & 4.97 \\
6053 & 377 & 12 & 4.85 & 6781 & 403 & 8 & 4.90 & 7541 & 428 & 6 & 4.93 & 8287
& 450 & 5 & 4.95 \\
6067 & 378 & 11 & 4.86 & 6791 & 403 & 9 & 4.90 & 7547 & 428 & 6 & 4.93 & 8291
& 453 & 2 & 4.98 \\
6073 & 377 & 12 & 4.84 & 6793 & 403 & 9 & 4.89 & 7549 & 428 & 6 & 4.93 & 8293
& 453 & 2 & 4.98 \\
6079 & 378 & 11 & 4.85 & 6803 & 402 & 10 & 4.88 & 7559 & 429 & 5 & 4.94 &
8297 & 453 & 2 & 4.98 \\
6089 & 378 & 12 & 4.85 & 6823 & 404 & 9 & 4.90 & 7561 & 429 & 5 & 4.94 & 8311
& 453 & 2 & 4.97 \\
6091 & 379 & 11 & 4.86 & 6827 & 404 & 9 & 4.89 & 7573 & 429 & 6 & 4.93 & 8317
& 454 & 1 & 4.98 \\
6101 & 379 & 11 & 4.86 & 6829 & 404 & 9 & 4.89 & 7577 & 428 & 7 & 4.92 & 8329
& 454 & 2 & 4.98 \\
6113 & 379 & 11 & 4.85 & 6833 & 404 & 9 & 4.89 & 7583 & 429 & 6 & 4.93 & 8353
& 454 & 2 & 4.97 \\ \hline
\end{array}
$\newpage

\textbf{Table 5. }The smallest known sizes
$\overline{t}_{2}=\overline{t} _{2}(2,q)$ of complete arcs in
planes $PG(2,q),$ $8363\leq q\leq 9109$, $ A_{q}=\left\lfloor
5\sqrt{q}-\overline{t}_{2}(2,q)\right\rfloor $, $ B_{q}\geq
\overline{t}_{2}(2,q)/\sqrt{q}$\smallskip

$\renewcommand{\arraystretch}{0.9}
\begin{array}{@{}r@{\,\,\,\,}c@{\,\,\,\,}c@{\,\,\,}c|@{\,\,}c@{\,\,\,\,}c@{\,\,\,\,}c@{\,\,\,}c|@{\,\,}c@{\,\,\,\,}c@{\,\,\,\,}c@{\,\,\,}c|@{\,\,}c@{\,\,\,\,}c@{\,\,\,\,}c@{\,\,\,\,}c}
\hline
q & \overline{t}_{2} & A_{q} & B_{q} & q & \overline{t}_{2} & A_{q} & B_{q}
& q & \overline{t}_{2} & A_{q} & B_{q} & q & \overline{t}_{2} & A_{q} &
B_{q}^{\phantom{H^{L}}} \\ \hline
8363 & 455 & 2 & 4.98 & 8573 & 462 & 0 & 4.99 & 8737 & 466 & 1 & 4.99 & 8929
& 471 & 1 & 4.99 \\
8369 & 456 & 1 & 4.99 & 8581 & 461 & 2 & 4.98 & 8741 & 467 & 0 & 5.00 & 8933
& 472 & 0 & 5.00 \\
8377 & 454 & 3 & 4.97 & 8597 & 463 & 0 & 5.00 & 8747 & 466 & 1 & 4.99 & 8941
& 472 & 0 & 5.00 \\
8387 & 456 & 1 & 4.98 & 8599 & 462 & 1 & 4.99 & 8753 & 467 & 0 & 5.00 & 8951
& 473 & 0 & 5.00 \\
8389 & 456 & 1 & 4.98 & 8609 & 463 & 0 & 5.00 & 8761 & 467 & 1 & 4.99 & 8963
& 473 & 0 & 5.00 \\
8419 & 457 & 1 & 4.99 & 8623 & 463 & 1 & 4.99 & 8779 & 468 & 0 & 5.00 & 8969
& 473 & 0 & 5.00 \\
8423 & 457 & 1 & 4.98 & 8627 & 464 & 0 & 5.00 & 8783 & 468 & 0 & 5.00 & 8971
& 472 & 1 & 4.99 \\
8429 & 457 & 2 & 4.98 & 8629 & 464 & 0 & 5.00 & 8803 & 468 & 1 & 4.99 & 8999
& 473 & 1 & 4.99 \\
8431 & 457 & 2 & 4.98 & 8641 & 464 & 0 & 5.00 & 8807 & 468 & 1 & 4.99 & 9001
& 474 & 0 & 5.00 \\
8443 & 457 & 2 & 4.98 & 8647 & 463 & 1 & 4.98 & 8819 & 469 & 0 & 5.00 & 9007
& 474 & 0 & 5.00 \\
8447 & 457 & 2 & 4.98 & 8663 & 465 & 0 & 5.00 & 8821 & 469 & 0 & 5.00 & 9011
& 472 & 2 & 4.98 \\
8461 & 459 & 0 & 5.00 & 8669 & 465 & 0 & 5.00 & 8831 & 469 & 0 & 5.00 & 9013
& 474 & 0 & 5.00 \\
8467 & 458 & 2 & 4.98 & 8677 & 465 & 0 & 5.00 & 8837 & 470 & 0 & 5.00 & 9029
& 475 & 0 & 5.00 \\
8501 & 459 & 2 & 4.98 & 8681 & 464 & 1 & 4.99 & 8839 & 470 & 0 & 5.00 & 9041
& 475 & 0 & 5.00 \\
8513 & 460 & 1 & 4.99 & 8689 & 465 & 1 & 4.99 & 8849 & 470 & 0 & 5.00 & 9043
& 475 & 0 & 5.00 \\
8521 & 460 & 1 & 4.99 & 8693 & 465 & 1 & 4.99 & 8861 & 470 & 0 & 5.00 & 9049
& 475 & 0 & 5.00 \\
8527 & 460 & 1 & 4.99 & 8699 & 465 & 1 & 4.99 & 8863 & 470 & 0 & 5.00 & 9059
& 475 & 0 & 5.00 \\
8537 & 460 & 1 & 4.98 & 8707 & 466 & 0 & 5.00 & 8867 & 470 & 0 & 5.00 & 9067
& 476 & 0 & 5.00 \\
8539 & 460 & 2 & 4.98 & 8713 & 466 & 0 & 5.00 & 8887 & 469 & 2 & 4.98 & 9091
& 477 &  & 5.01 \\
8543 & 461 & 1 & 4.99 & 8719 & 466 & 0 & 5.00 & 8893 & 471 & 0 & 5.00 & 9103
& 479 &  & 5.03 \\
8563 & 461 & 1 & 4.99 & 8731 & 466 & 1 & 4.99 & 8923 & 472 & 0 & 5.00 & 9109
& 479 &  & 5.02 \\ \hline
\end{array}
$
\end{center}

\section{Observations on $\overline{t}_{2}(2,q)$ values}

\label{sec4_observations}

We look for upper estimates of the collection of
$\overline{t}_{2}(2,q)$ values from Tables 1-5 in the
form~(\ref{eq1_KimVu_c=300}), see \cite{KV},
\cite[Tab.\thinspace 2.6]{HirsStor-2001}, and
\cite[Sec.\thinspace 4]{BDFMP-DM}. For definiteness, we use the
natural logarithms. Let $c$ be a constant independent of $q$.
We introduce $ D_{q}(c)$ and $\overline{D}_{q}(c)$ as follows:
\begin{eqnarray}
t_{2}(2,q) &=&D_{q}(c)\sqrt{q}\ln ^{c}q,~  \notag \\
\overline{t}_{2}(2,q) &=&\overline{D}_{q}(c)\sqrt{q}\ln ^{c}q.
\label{eq4_Dq(c)}
\end{eqnarray}

Let $\overline{D}_{\text{aver}}(c,q_{0})$ be the average value of $\overline{
D}_{q}(c)$ calculated in the region $q_{0}\leq q\leq 9109$ under condition $
q\notin N.$

From Tables 1-5, we obtain Observation~1.\smallskip

\noindent \textbf{Observation 1.} \emph{Let }$173\leq q\leq 9109$ \emph{\
under condition }$q\notin N$\emph{. Then when $q$ grows}, \emph{the values
of }$\overline{D}_{q}(0.75)$\emph{\ oscillate about the average value }$
\overline{D}_{\text{aver}}(0.75,173)=0.95579$ (see
Fig.~1).\emph{\ Also,}
\begin{equation}
\begin{array}{cl}
0.946<\overline{D}_{q}(0.75)<0.9634 & \mbox{if }173\leq q<1000,\smallskip
\\
0.953<\overline{D}_{q}(0.75)<0.9605 & \mbox{if }1000<q<2000,\smallskip  \\
0.950<\overline{D}_{q}(0.75)<0.9595 & \mbox{if }2000<q<3000,\smallskip  \\
0.950<\overline{D}_{q}(0.75)<0.9588 & \mbox{if }3000<q<4000,\smallskip  \\
0.951<\overline{D}_{q}(0.75)<0.9584 & \mbox{if }4000<q<5000,\smallskip  \\
0.950<\overline{D}_{q}(0.75)<0.9579 & \mbox{if }5000<q<6000,\smallskip  \\
0.951<\overline{D}_{q}(0.75)<0.9577 & \mbox{if }6000<q<7000,\smallskip  \\
0.947<\overline{D}_{q}(0.75)<0.9573 & \mbox{if }7000<q<8000,\smallskip  \\
0.949<\overline{D}_{q}(0.75)<0.9573 & \mbox{if }8000<q.
\end{array}
\label{eq4_Dq74}
\end{equation}
\begin{figure}[t]
\epsfig{file=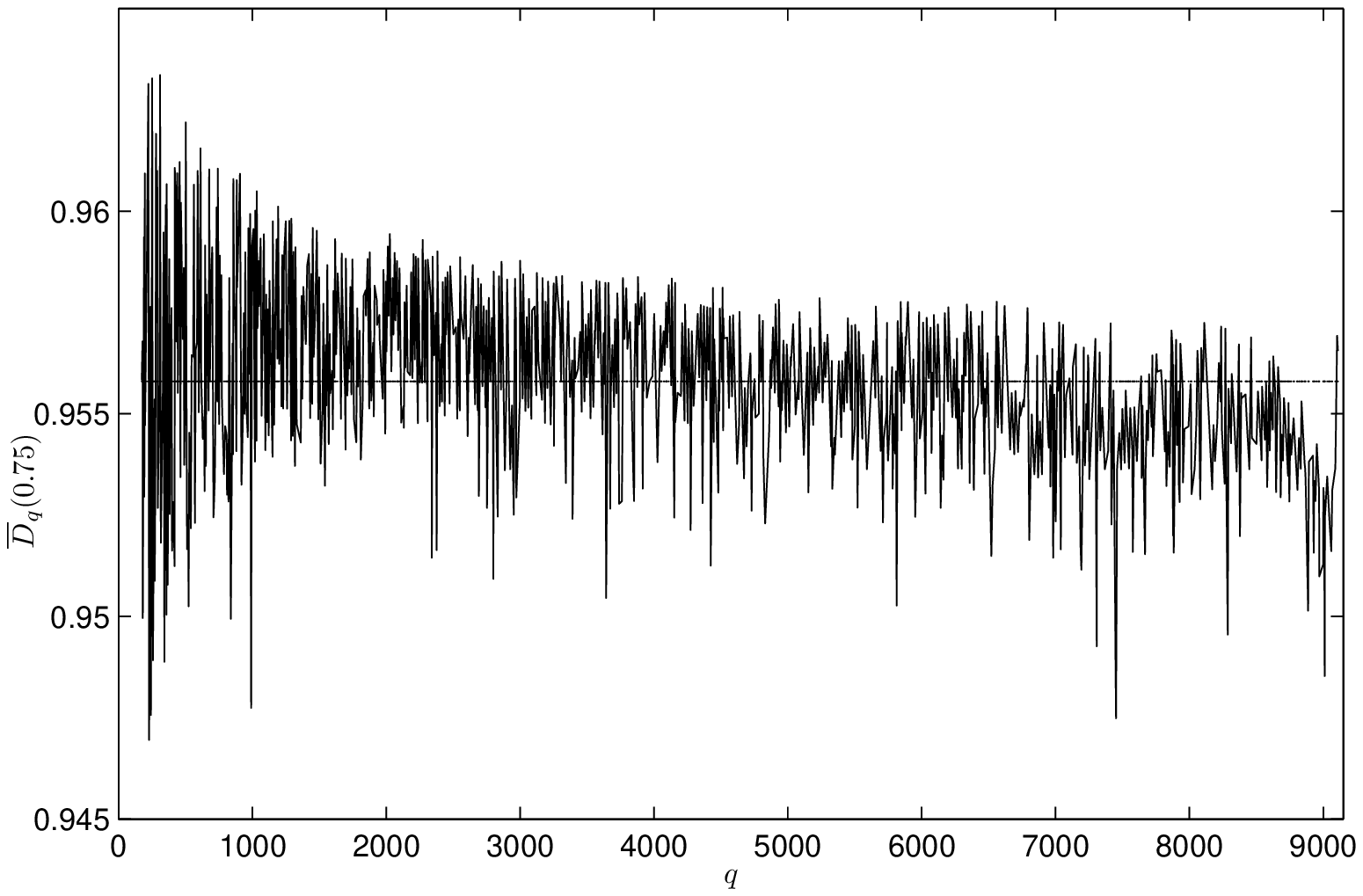,width=\textwidth}
\caption{The values of $\overline{D}_{q}(0.75)$ for $173\le q\le 9109$, $ q\notin N$.
$\overline{ D}_{\text{aver}}(0.75,173)=0.95579$}
\end{figure}
By Observation 1 it seems that the values of $D_{q}(0.75)$ and
$\overline{D} _{q}(0.75)$ are sufficiently convenient for
estimates of $t_{2}(2,q)$ and $ \overline{t}_{2}(2,q).$

From Tables 1-5, we obtain Theorem \ref{th4_ln0.75}.

\begin{theorem}
\label{th4_ln0.75} In $PG(2,q),$
\begin{equation}
t_{2}(2,q)<0.9987\sqrt{q}\ln ^{0.75}q\quad \mbox{ for }23\leq q\leq 9109,
\text{ }q\in T_{3}.  \label{eq4_9987}
\end{equation}
\end{theorem}

The graphs of values of $\sqrt{q}\ln ^{0.8}q$, $\sqrt{q}\ln
^{0.75}q$, $ \overline{t}_{2}(2,q)$, and $\sqrt{q}\ln ^{0.5}q$
are shown \smallskip on Fig.~2 where $\sqrt{q}\ln ^{0.8}q$ is
the top curve and $\sqrt{q}\ln ^{0.5}q$
\smallskip is the bottom one.
\begin{figure}[t]
\epsfig{file=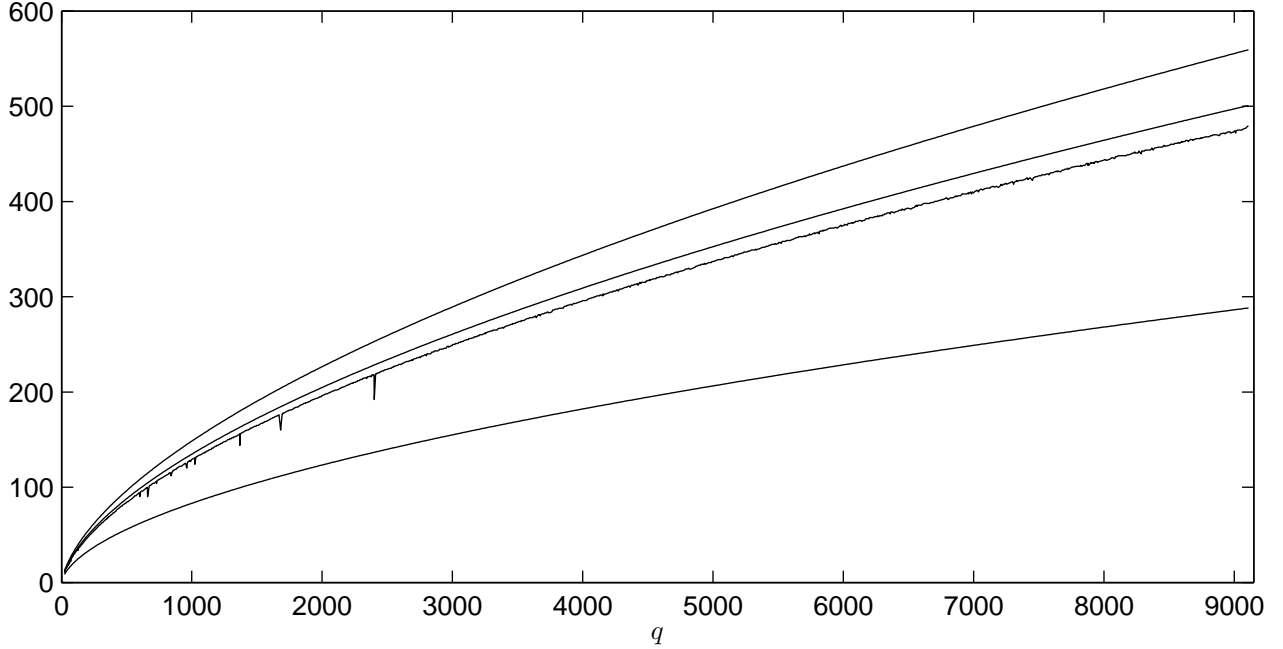,width=\textwidth}
\caption{The values of $\sqrt{q}\ln^{0.8}q$ (the top curve),
$\sqrt{q}\ln^{0.75}q$ (the 2-nd curve), $\overline{t}_{2}(2,q)$ (the 3-rd curve), and
$\sqrt{q}\ln^{0.5}q$ (the bottom curve) for $23\le q\le 9109$}
\end{figure}

One can see on Fig.~2 that always $\overline{t}_{2}(2,q)<\sqrt{q}\ln^{0.75}q$
\smallskip and, moreover, when $q$ grows, the graphs $\sqrt{q}\ln^{0.75}q$
and $\overline{t}_{2}(2,q)$ diverge so that positive difference $\sqrt{q}
\ln^{0.75}q - \overline{t}_{2}(2,q)$ increases.

We denote
\begin{equation}
\widehat{t}_{2}(2,q)=\overline{D}_{\text{aver}}(0.75,173)\sqrt{q}\ln
^{0.75}q,\quad \overline{\Delta }_{q}=\overline{t}_{2}(2,q)-\widehat{t}
_{2}(2,q),\quad \overline{P}_{q}=\frac{100\overline{\Delta }_{q}}{\overline{
t }_{2}(2,q)}\%.  \label{eq4_Delta-q_Pq}
\end{equation}
One can treat $\widehat{t}_{2}(2,q)$ as a \emph{predicted}
value of $ t_{2}(2,q)$. Then $\overline{\Delta }_{q}$ is the
difference between the smallest known size
$\overline{t}_{2}(2,q)$ of complete arcs and the predicted
value. Finally, $\overline{P}_{q}$ is this difference in
percentage terms of the smallest known size.\smallskip

\noindent \textbf{Observation 2.} \emph{Let }$173\leq q\leq 9109$, $q\notin N
$\emph{. Then\newline
}
\begin{equation}
-3.70<\overline{\Delta }_{q}<0.81.  \label{eq4_Delta_q}
\end{equation}

\begin{equation}
\begin{array}{cl}
-0.94\%<\overline{P}_{q}<0.79\% & \mbox{if }173\leq q<1000,\smallskip  \\
-0.28\%<\overline{P}_{q}<0.49\% & \mbox{if }1000<q<2000,\smallskip  \\
-0.52\%<\overline{P}_{q}<0.38\% & \mbox{if }2000<q<3000,\smallskip  \\
-0.57\%<\overline{P}_{q}<0.32\% & \mbox{if }3000<q<4000,\smallskip  \\
-0.48\%<\overline{P}_{q}<0.27\% & \mbox{if }4000<q<5000,\smallskip  \\
-0.59\%<\overline{P}_{q}<0.22\% & \mbox{if }5000<q<6000,\smallskip  \\
-0.46\%<\overline{P}_{q}<0.20\% & \mbox{if }6000<q<7000,\smallskip  \\
-0.88\%<\overline{P}_{q}<0.16\% & \mbox{if }7000<q<8000,\smallskip  \\
-0.66\%<\overline{P}_{q}<0.16\% & \mbox{if }8000<q.
\end{array}
\label{eq4_percent}
\end{equation}

By (\ref{eq4_Delta_q}) and (\ref{eq4_percent}), see also Fig. 3
and 4, the upper bounds of $\overline{\Delta }_{q}$ and
$\overline{P}_{q}$ are relatively small. Moreover, the upper
bound of $\overline{P}_{q}$ decreases when $q$ grows. Therefore
the values of $\overline{\Delta }_{q}$ and $ \overline{P}_{q}$
are useful for computer search of small arcs.

\begin{figure}[tb]
\epsfig{file=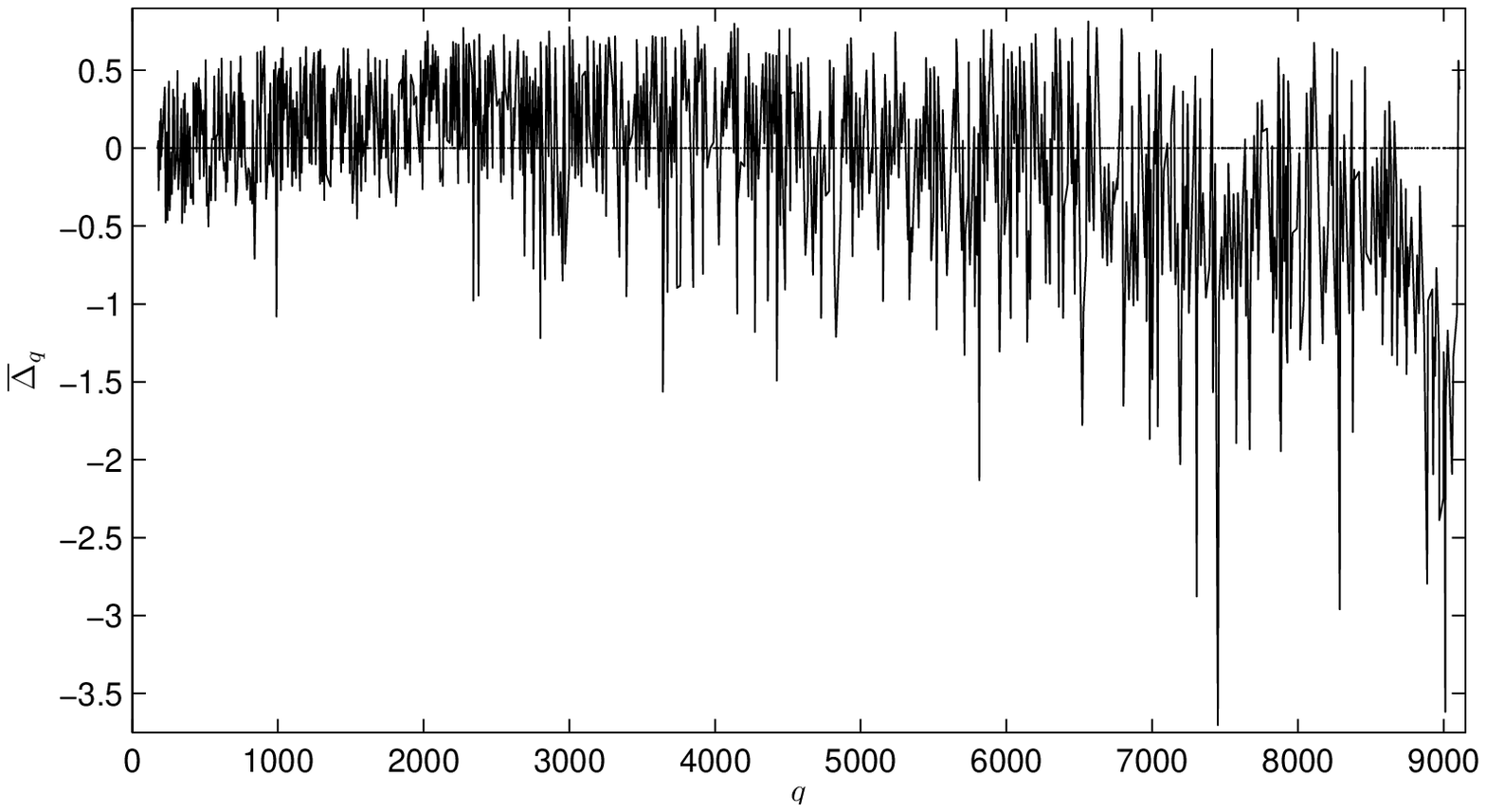,width=\textwidth}
\caption{The values of $\overline{\Delta }_{q}=\overline{t}_{2}(2,q)-\widehat{t}
_{2}(2,q)$ for $173\le q\le 9109$,  $ q\notin N$}
\end{figure}
\begin{figure}[tb]
\epsfig{file=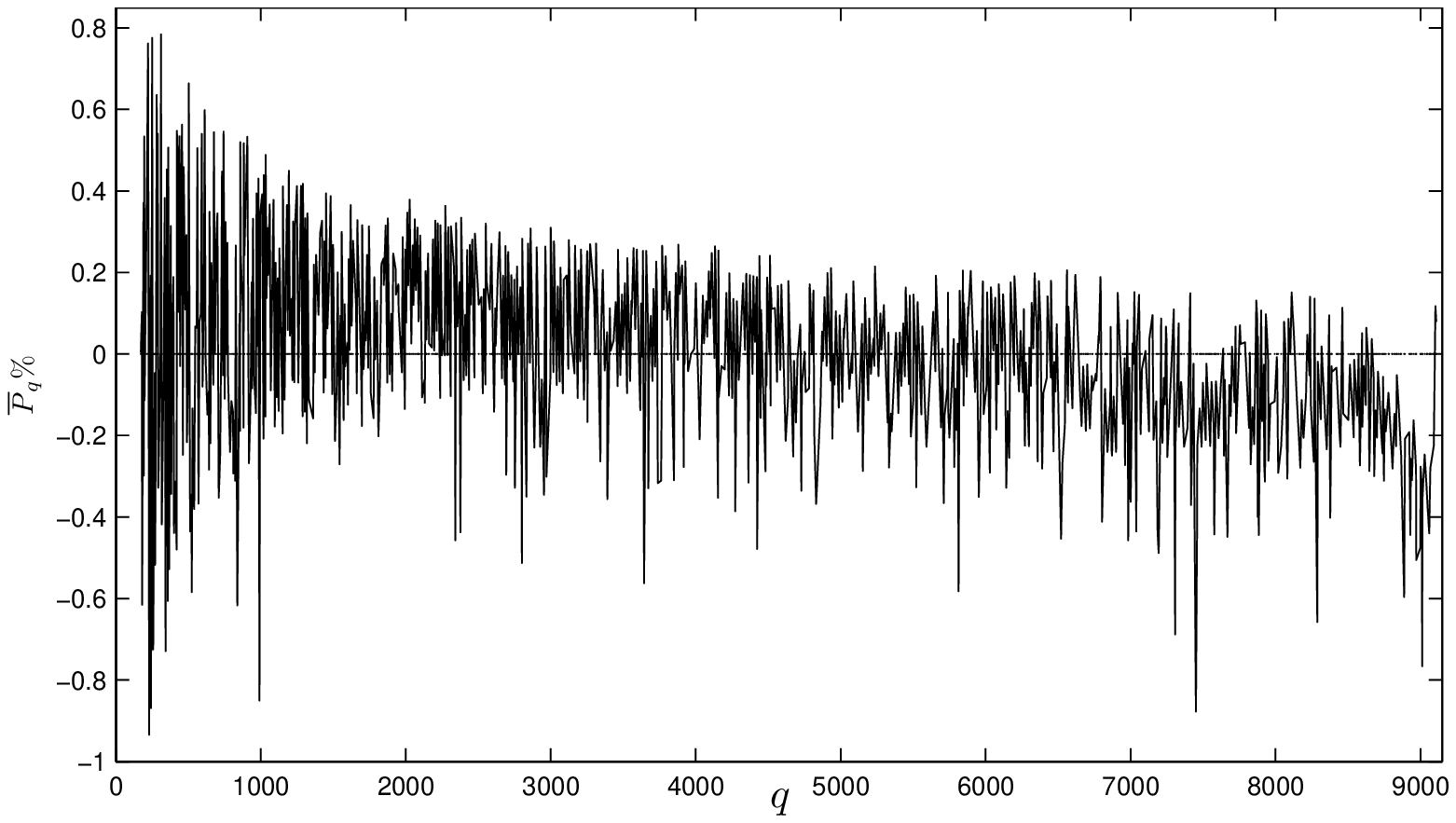,width=\textwidth}
\caption{The values of $\overline{P}_{q}=100\overline{\Delta }_{q}/\overline{
t }_{2}(2,q)\%$ for $173\le q\le 9109$,  $ q\notin N$}
\end{figure}
The relations (\ref{eq4_Dq74})--(\ref{eq4_percent}), Theorems
\ref {th1_ln0.75} and \ref{th4_ln0.75}, and Figures 1--4 are
the foundation for Conjecture \ref{conj1_ ln0.75}.

\begin{remark}
\label{rem4} By above, $\sqrt{q}\ln ^{0.75}q$\, seems to be a
reasonable upper bound on the current collection of
$\overline{t}_{2}(2,q)$ values. It gives some reference points
for computer search and foundations for Conjecture \ref {conj1_
ln0.75} on the upper bound for $t_{2}(2,q)$. In principle, the
constant $c=0.75$ can be sightly reduced to move the curve
$\sqrt{q}\ln ^{c}q $ near to the curve of
$\overline{t}_{2}(2,q)$, see Fig.~2 and \cite[ Th.\thinspace
4.3]{BDFMP-DM}.
\end{remark}

\end{document}